\documentclass {article}
\usepackage{amssymb}
\usepackage{amsmath}
\usepackage{graphicx}
\usepackage{color}
\newtheorem{theorem}{Theorem}

\newtheorem{definition}{Definition}
\newtheorem{lemma}[theorem]{Lemma}
\newtheorem{rem}[theorem]{Remark}

\newcommand{\inte }{{\rm int}\,}

\newcommand{\cl }{{\rm cl}\,}

\newcommand{\cover}[1]{\stackrel{#1}{\Longrightarrow}}
\newcommand{\comment}[1]{\mbox{}}

\def\qed{{\hfill{\vrule height5pt width3pt depth0pt}\medskip}}

\begin{document}
\begin{center} {\bf \LARGE Topological shadowing and  the Grobman-Hartman Theorem} \\
\vskip 0.5cm {\large Piotr Zgliczy\'nski}\footnote{Research has
been supported by Polish
National Science Centre grant 2011/03B/ST1/04780} \\
Jagiellonian University, Institute of Computer Science and Computational Mathematics, \\
{\L}ojasiewicza 6, 30--348  Krak\'ow, Poland \\

 \vskip 0.5cm
e-mail:  piotr.zgliczynski@ii.uj.edu.pl

 \vskip 0.5cm
 \today
\end{center}

\begin{abstract}
We give geometric proofs for Grobman-Hartman theorem for
diffeomorphisms and ODEs. Proofs use covering relations and cone
conditions for maps and isolating segments and cone conditions for
ODEs. We establish a topological versions of the Grobman-Hartman theorem as the existence of some
semiconjugaces.
\end{abstract}



\section{Introduction}

The goal of this paper is to give a new geometric proof of the
Grobman-Hartman \cite{G1,G2,H1} theorem for diffeomorphism and ODEs in finite
dimension. By 'the geometric proof' we understand  the proof
which works in the phasespace of the system under consideration
and uses concepts of qualitative geometric nature.

 We focus on the global version of the Grobman-Hartman theorem, which  in the case map states that, if $A:\mathbb{R}^n \to \mathbb{R}^n$ is a hyperbolic linear isomorphism and if $g:\mathbb{R}^n \to \mathbb{R}^n$ is given by
\begin{equation}
  g(x)=Ax + h(x),
\end{equation}
where $h:\mathbb{R}^n \to \mathbb{R}^n$ is a bounded $C^1$ function, such that $\|Dh(x)\| \leq \epsilon$ for $x \in \mathbb{R}^n$, then if $\epsilon$ is sufficiently small, then $A$ and $g$ are conjugated by a continuous homeomorphism.

There are many of proofs of the Grobman-Hartman theorem in the literature. An exemplary geometric proof
 can be found in the Katok-Hasselblatt book \cite{KH}. This proof is placed in the context of the hyperbolicity, they show that
 dynamics of $g$  is hyperbolic on whole $\mathbb{R}^n$ and the conjugating homeomorphism is constructed geometrically by considering the stable and unstable
 leaves of  points to construct the linearizing coordinates.

The other family of proofs of the Grobman-Hartman theorem uses tools from the functional analysis.
 The standard functional analysis proof \cite{Pa,Pu,BV}, which
is now a textbook proof (see for example \cite{A,C99,PdM,Ze}), it
studies the conjugacy problem in some abstract Banach space of
maps. The original proof by P. Hartman \cite{H1,H2,H3} also
belongs to this category,  but it lacks  the simplicity of the contemporary approach, because to solve the
conjugacy problem Hartman required first to introduce new
coordinates which straighten the invariant manifolds of the hyperbolic fixed point. The standard functional analysis proof,
whose idea apparently comes from paper by Moser \cite{M}(see also \cite{Pa,Pu}), in a
current form is a straightforward application of the Banach contraction principle. The whole effort is to chose the correct
Banach space and a contraction, whose fixed point will give us the
conjugacy.


 In this paper we would like to give a new geometric proof  the global version of the Grobman-Hartman theorem (Theorem~\ref{thm:GH-maps}). The geometric idea behind our
approach can be seen a shadowing of $\delta$-pseudo orbit, with $\delta$ not small. This is accomplished using covering
relations and the cone condition~\cite{ZGi,ZCC} in case of
diffeomorphisms and for ODEs the notion of the isolating
segment~\cite{S1,S2,S3,SW,WZ} and the cone conditions has been
used. Compared to the geometric proof in \cite{KH} we stress more the topological aspects.  As the byproduct of our approach we obtain
two topological variants of the Grobman-Hartman theorem
\begin{itemize}
\item if we drop the assumption that $\|Dh\|$ is small, but we demand instead that $g$ is homeomorphism, then we show that there exists a semiconjugacy
between $A$ and $g$, see Theorem~\ref{thm:GH-maps-top} for the precise statement,
\item if we drop the assumption that $\|Dh\|$ is small, then we show  that there exists a semiconjugacy
between $A$ restricted to the unstable subspace and $g$, see Theorem~\ref{thm:GH-maps-top-no-inj} for the precise statement,
\end{itemize}

 Let us comment about the relation between our proofs of the
theorem  for maps and for ODEs. The standard approach would be to
derive the ODE case from the map case, by considering the time
shift by one time unit and then arguing that  we can obtain from
it the conjugacy for all times (see \cite{H1,Pa,Pu,PdM}). Here, we
provide the proof for ODEs which is independent from the map case
in order to illustrate the power of the concept the isolating
segment with the aim to obtain a clean ODE-type proof. For an another clean ODE-type proof
using the functional analysis type arguments see \cite{CS}.


Regarding the regularity of the conjugating homeomorphism in the global Grobman-Hartman theorem  there is a nice  argument of geometric nature in Katok-Hasselblatt book \cite{KH} that shows that this conjugacy between  has to be H\"older.  However  no effort is made there
to estimate the H\"older exponent. Using our shadowing ideas we estimate  this exponent. We
obtain the same estimate for the H\"older exponent as in
the work by Barreira and Valls \cite{BV}, Belitskii \cite{B}, Belitskii and Rayskin \cite{BR} which  apparently are the best results in this directions (see \cite{BV} and references given there). In these papers  the functional analysis type of reasoning was used and results are valid also in the Banach space.

The organization of this paper can be described as follows.
Section~\ref{sec:GHmapsglobal} contains the geometric proof of the
global version of the Grobman-Hartman theorem. In
Section~\ref{sec:GH-Holder} we show the H\"older regularity of the
conjugacy in the Grobman-Hartman theorem.

Section~\ref{sec:ODE-GH} contains a geometric proof of the
Grobman-Hartman theorem for flows, which is independent from the
proof for maps.

At the end of this paper we included two appendices, which
contains relevant definitions and theorems about the covering
relations and the isolating segments.

\subsection{Notation}
If $A \in \mathbb{R}^{d_1 \times d_2}$ is a matrix, then by $A^t$ we
will denote its transpose. By $B(x,r)$ we will denote the open
ball centered at $x$ and radius $r$. For maps depending on some
parameters $h:P \times X \to X$ by $h_p:X \to X$ we will denote
the map $h_p(x)=h(p,x)$.

In this note we will work in $\mathbb{R}^n=\mathbb{R}^u \times
\mathbb{R}^s$. According to  this decomposition we will often
represent points $z \in \mathbb{R}^n$ as $z=(x,y)$, where $x \in
\mathbb{R}^u $ and $y \in \mathbb{R}^s$. On $\mathbb{R}^n$ we
assume the standard scalar product $(u,v)=\sum_i u_i v_i$. This
scalar product induces the norm on $\mathbb{R}^u$ and
$\mathbb{R}^s$. We will use the following norm on $\mathbb{R}^n$,
$\|(x,y)\|_{max}=\max (\|x\|,\|y\|)$ and we will usually drop the
subscript $max$.

We will use also projections $\pi_x$ and $\pi_y$, so that $\pi_x(x,y)=x$ and
$\pi_y(x,y)=y$.

\section{Global version of the Grobman-Hartman theorem for maps}
\label{sec:GHmapsglobal}

In this section we will give a geometric proof of the Grobman-Hartman
theorem for maps and its topological variants.

We will consider a map $g:\mathbb{R}^n \to \mathbb{R}^n$, such that
\begin{equation}
  g(z)=A(z)+h(z).
\end{equation}

We will have the following set of assumptions on $A$ and $h$, which we will refer to as the \emph{standard conditions}
\begin{itemize}
\item We assume that $A: \mathbb{R}^n \to \mathbb{R}^n$
is a linear isomorphism, of the following form
\begin{equation}
  A(x,y)=(A_u x, A_s y),
\end{equation}
where $n=u+s$, $A_u: \mathbb{R}^u \to \mathbb{R}^u$ and
$A_s:\mathbb{R}^s \to \mathbb{R}^s$ are linear isomorphisms such
that
\begin{eqnarray}
  \|A_u x\| &\geq& c_u \|x\|, \quad c_u >1, \quad \forall x \in \mathbb{R}^u \\
  \|A_s y \| &\leq& c_s \|y\|, \quad 0 < c_s < 1, \quad \forall y \in
  \mathbb{R}^s.
\end{eqnarray}
\item map $h:\mathbb{R}^n \to \mathbb{R}^n$ is continuous and  there exist $M$ such that
\begin{eqnarray}
  \|h(x)\| \leq M, \quad \forall x \in \mathbb{R}^n
\end{eqnarray}
\end{itemize}

\begin{theorem}
\label{thm:GH-maps}
Assume the standard conditions.

Additionally assume that $h$ is of class $C^1$ and
such that there exist  $\epsilon$ such that
\begin{eqnarray}
  \|Dh(x)\| \leq \epsilon, \quad \forall x \in \mathbb{R}^n.
\end{eqnarray}

Then there exists
$\epsilon_0=\epsilon_0(A)>0$, such that if $\epsilon <
\epsilon_0(A)$, then there exists a homeomorphism
$\sigma:\mathbb{R}^n \to \mathbb{R}^n$ such that
\begin{equation}
    \sigma \circ g = A \circ \sigma.  \label{eq:GH-conj}
\end{equation}

\end{theorem}

\textbf{Comment:} Observe that there is no bound on $M$, we also do not assume that $h(0)=0$.

In the next theorem we drop the assumption that $h$ is $C^1$ with small $Dh$, but we keep the requirement that $g$ is an injective map.

\begin{theorem}
\label{thm:GH-maps-top}
Assume the standard conditions.

  Assume  map $g$ is an injection.

Then there exists a continuous surjective  map
$\sigma:\mathbb{R}^n \to \mathbb{R}^n$ such that
\begin{equation}
   \sigma \circ g = A \circ \sigma.  \label{eq:GH-conj-top}
\end{equation}

\end{theorem}

In the next theorem  we will drop the assumption that $g$ is an injection. Then we we no longer have
a unique full trajectory through a point for map $h$.

\begin{theorem}
\label{thm:GH-maps-top-no-inj} Assume the standard conditions.

Then there exists a continuous surjective  map
$\sigma_u:\mathbb{R}^n \to \mathbb{R}^u$ such that
\begin{equation}
   \sigma_u \circ g = A_u \circ \sigma_u.  \label{eq:GH-conj-top-no-inj}
\end{equation}

\end{theorem}

Before the proof of Theorem~\ref{thm:GH-maps}, \ref{thm:GH-maps-top}, \ref{thm:GH-maps-top-no-inj} we need first to
develop some technical tools. The basic steps and constructions
used in the proof are given in
Section~\ref{subsec:thm-gh-maps-proof}. We invite the reader to
jump first to this section to see the overall picture of the proof
and then consult other more technical sections when necessary.

We will use the following notation: $g_\lambda=A+\lambda h$ for
$\lambda \in [0,1]$. In this notation  we have $g=g_1$.

\subsection{$g_\lambda$ are onto }
\begin{lemma}
\label{lem:onto}
Assume standard conditions. Then $g_\lambda$ are onto, i.e $g_\lambda(\mathbb{R}^n)=\mathbb{R}^n$.
\end{lemma}
\textbf{Proof:}
 The surjectivity of $g_\lambda$ follows from the
following observation: a bounded continuous perturbation a linear
isomorphism is a surjection - the proof is based on the local
Brouwer degree (see for example Appendix in \cite{ZGi} for the
definition and properties). Details are as follows.

For fixed $y \in \mathbb{R}^n$ we consider equation
$y=g_\lambda(x)$, which is equivalent to  $x+\lambda
A^{-1}h(x)=A^{-1}y=\tilde{y}$. Let us define a map
\begin{equation}
  F_\lambda(x)=x + \lambda A^{-1}h(x) - \tilde{y}.
\end{equation}
Observe that if $\|x - \tilde{y}\|  > \|A^{-1}\| M$, then
$F_\lambda(x) \neq 0$.

This shows that
$\deg(F_\lambda,B(\tilde{y,}\|A^{-1}\| M ),0)$  (the local Brouwer
degree of $F_\lambda$ on the set $B(\tilde{y,}\|A^{-1}\| M )$ at
$0$ ) is defined and
\begin{equation}
 \deg(F_\lambda,B(\tilde{y,}\|A^{-1}\| M ),0) = \deg(F_0,B(\tilde{y,}\|A^{-1}\| M
 ),0), \quad \forall \lambda \in [0,1].
\end{equation}
But for $\lambda=0$ we have $F_0(x)=x - \tilde{y}$. Hence
$\deg(F_0,B(\tilde{y,}\|A^{-1}\|M,0)=1$. Therefore
$F_\lambda(x)=0$ has solution for any $\tilde{y} \in
\mathbb{R}^n$.
\qed

\subsection{$g_\lambda$ are homeomorphisms under assumptions of Theorem~\ref{thm:GH-maps}}

The following lemma can be found for example in  \cite[Lemma 1]{Pu} \cite[Proposition II.2]{Ze}
\begin{lemma}
\label{lem:homeo} Let $A$ and $h$ be as in Theorem~\ref{thm:GH-maps}. Let $\epsilon_1(A)=\frac{1}{\|A^{-1}\|}
>0$.

 If $\epsilon < \epsilon_1(A)$, then $g_\lambda$
is a homeomorphism and $g_\lambda^{-1}$ is Lipschitz.
\end{lemma}
\noindent \textbf{Proof:}
The surjectivity follows from Lemma~\ref{lem:onto}.

The injectivity is obtained as follows
\begin{eqnarray*}
  \|g_\lambda(z_1) - g_\lambda(z_2)\| = \|Az_1 + \lambda h(z_1) - (Az_2 + \lambda
  h(z_2))\| \geq \\
   \|A(z_1) - A(z_2) \| - \lambda \|h(z_1) -
  h(z_2)\| \geq \\
  \frac{1}{\|A^{-1}\|} \|z_1-z_2\| - \epsilon
  \|z_1-z_2\|= \left(\frac{1}{\|A^{-1}\|}  - \epsilon\right)
  \|z_1-z_2\|.
\end{eqnarray*}

From the above formula it follows also that
\begin{equation}
 \|z_1 - z_2\| \geq \left(\frac{1}{\|A^{-1}\|}  - \epsilon\right)
  \|g_\lambda^{-1}(z_1)-g_\lambda^{-1}(z_2)\|.
\end{equation}
Therefore
\begin{equation}
  \|g_\lambda^{-1}(z_1)-g_\lambda^{-1}(z_2)\| \leq \left(\frac{1}{\|A^{-1}\|}  - \epsilon\right)^{-1} \|z_1 - z_2\|
\end{equation}

 \qed

\subsection{Cone condition for $g_\lambda$ under assumptions of Theorem~\ref{thm:GH-maps}}
\label{subsec:cc-maps}
Throughout this subsection we work under assumptions of Theorem~\ref{thm:GH-maps}.

 We will establish the cone
condition for $g_\lambda$ using the approach from~\cite{ZCC},
where the cones are defined in terms of a quadratic form.

 Let $Q$ be an quadratic form in $\mathbb{R}^n=\mathbb{R}^u \times \mathbb{R}^s$ given by $Q(x,y)=(x,x)  - (y,y)$. Our goal is to show the following
\emph{cone condition}: for sufficiently small $\eta
>0$ it holds
\begin{equation}
  Q(A z_1 - Az_2) > (1\pm \eta) Q(z_1 - z_2), \quad z_1,z_2 \in
  \mathbb{R}^n, z_1 \neq z_2.  \label{eq:GH-ccA}
\end{equation}
This will be established in Lemma~\ref{lem:cc}.

By $Q$  we will also denote a matrix, such that $Q(z)=z^tQz$. In
our case $Q=\left[\begin{array}{cc}
  I_u & 0 \\
  0 & -I_s
\end{array}
\right]$,  where $I_u \in \mathbb{R}^{u \times u}$ and $I_s \in
\mathbb{R}^{s \times s}$ are the identity matrices.

\begin{lemma}
\label{lem:ccA}
 For $0\leq \eta \leq \min(c_u^2 -1, 1-c_s^2)$ the matrix $A^tQA - (1\pm \eta)Q$ is positive definite.
\end{lemma}
\textbf{Proof:} Easy computations show that
\begin{eqnarray*}
  A^tQA = \begin{pmatrix}
    A_u^tA_u & 0 \\
    0 & A_s^t A_s \
  \end{pmatrix}.
\end{eqnarray*}
Hence for any $z=(x,y) \in \mathbb{R}^u \times \mathbb{R}^s
\setminus \{0\}$ holds
\begin{eqnarray*}
   z^t\left(A^tQA - (1\pm \eta)Q\right) z = x^tA_u^tA_ux - (1\pm \eta)x^2 +
   (1\pm\eta)y^2 - y^tA_s^tA_sy = \\
   (A_ux,A_ux) -  (1\pm \eta)x^2 + (1\pm\eta)y^2 - (A_sy,A_sy)
   \geq \\
    (c_u^2 - 1 - \eta)x^2 + (1-\eta - c_s^2) y^2 >0,
\end{eqnarray*}
if $c_u^2-1 > \eta$ and $1-c_s^2 > \eta$.
 \qed

\begin{lemma}
\label{lem:cc} There exists $\epsilon_0(A)>0$, such that if $0\leq
\epsilon < \epsilon_0(A)$, then  there exists $\eta \in (0,1)$
such that for any $\lambda \in [0,1]$ the following \emph{cone
condition} holds
\begin{equation}
  Q(g_\lambda(z_1) - g_\lambda(z_2)) > (1\pm \eta) Q(z_1-z_2),
  \quad \forall z_1,z_2\in \mathbb{R}^n, z_1 \neq z_2.  \label{eq:cc}
\end{equation}
\end{lemma}
\textbf{Proof:}
 We have
\begin{eqnarray*}
  Q(g_\lambda(z_1) - g_\lambda(z_2)) = (z_1-z_2)^t (D(z_1,z_2)^t Q D(z_1,z_2))
  (z_1-z_2), \\
  D(z_1,z_2)=  \int_0^1 Dg_\lambda(t(z_1-z_2) + z_2)dt
\end{eqnarray*}
Let
\begin{equation}
 C(z_1,z_2)=  \int_0^1 Dh(t(z_1-z_2) + z_2)dt,
\end{equation}
then
\begin{equation}
  D(z_1,z_2)=A+ \lambda C(z_1,z_2).
\end{equation}

 Observe that $\|C(z_1,z_2)\| \leq \epsilon$.

From Lemma~\ref{lem:ccA} it follows that $A^t Q A - (1\pm \eta) Q$
is positive definite for sufficiently small $\eta >0$. Let us fix
such $\eta$.

Since being a positively defined symmetric matrix is an open
condition, hence there exists $\epsilon_0(A)>0$ be such that the
matrix
\begin{equation}
  (A+\lambda C)^t Q (A+\lambda C) - (1\pm \eta) Q  \label{eq:Acc}
\end{equation}
is positive definite for any $\lambda \in [0,1]$ and $C \in
\mathbb{R}^{n \times n}$ satisfying $\|C\| \leq \epsilon_0$. \qed

From Lemma~\ref{lem:homeo} it follows that for any $\lambda \in
[0,1]$ and any point $z$ we can define a full orbit for
$g_\lambda$ through this point, i.e. $g_\lambda^{k}(z)$ makes
sense for any $k \in \mathbb{Z}$.

\begin{lemma}
\label{lem:bdiff-uniqueness} Assume that $\epsilon <
\min(\epsilon_0(A),\epsilon_1(A))$ from Lemmas~\ref{lem:cc} and
\ref{lem:homeo}. Let $\lambda \in [0,1]$. If $z_1,z_2 \in
\mathbb{R}^n$ and $\beta$ are such that
\begin{equation}
  \|g^k_\lambda(z_1) - g^k_{\lambda}(z_2)\| \leq \beta, \quad
  \forall k \in \mathbb{Z}, \label{eq:bounded-diff}
\end{equation}
then $z_1=z_2$.
\end{lemma}
\textbf{Proof:} The proof is by  contradiction. Assume that
$z_1 \neq z_2$. Either $Q(z_1 - z_2) \geq 0$ or $Q(z_1 - z_2) <0$.

Let us consider first  case $Q(z_1 - z_2) \geq 0$. By the cone
condition (Lemma~\ref{lem:cc}) we obtain for any $k>0$
\begin{eqnarray*}
  Q(g_\lambda(z_1) - g_\lambda(z_2)) &>& Q(z_1 - z_2) \geq 0 \\
 \|\pi_x (g_\lambda^k(z_1) - g_\lambda^k(z_2)) \| &\geq& Q(g_\lambda^k(z_1) - g_\lambda^k(z_2)) >
  (1+\eta)^{k-1}Q(g_\lambda(z_1) - g_\lambda(z_2)).
\end{eqnarray*}
Therefore $g_\lambda^k(z_1) - g_\lambda^k(z_2)$ is unbounded. This
contradicts \eqref{eq:bounded-diff}.

Now we consider the case $ Q(z_1 - z_2) <0$. The cone condition (Lemma~\ref{lem:cc})
applied to the inverse map gives for any $k >0$
\begin{eqnarray*}
  Q(z_1 - z_2) > (1-\eta) Q(g_\lambda^{-1}(z_1) - g_\lambda^{-1}(z_2)) > \\
    (1-\eta)^k Q(g_\lambda^{-k}(z_1) - g_\lambda^{-k}(z_2)).
\end{eqnarray*}
Therefore we obtain
\begin{equation}
  -Q(g_\lambda^{-k}(z_1) - g_\lambda^{-k}(z_2)) > \frac{1}{(1-\eta)^k}(-Q(z_1 -
  z_2)).
\end{equation}
Therefore $g_\lambda^{-k}(z_1) - g_\lambda^{-k}(z_2)$ is
unbounded. This contradicts \eqref{eq:bounded-diff}. \qed

\subsection{Covering relations}

We assume that the reader is familiar with the notion of h-set and
covering relation \cite{ZGi}. For the convenience of the reader we
recall these notions in Appendix~\ref{app:cov-rel}.

\begin{definition}
\label{def:N(a,b)}
For any $z \in \mathbb{R}^n$, $\alpha >0$ we define an h-set (with
a natural structure)  $N(z,\alpha)=z + \overline{B}_u(0,\alpha)
\times \overline{B}_s(0,\alpha)$.
\end{definition}

The following theorem follows immediately from
Theorem~\ref{thm:inf-chain} in Appendix \ref{app:cov-rel}.
\begin{theorem}
\label{thm:cov-rel-exists}  Assume that we have a bi-infinite chain
of covering relations
\begin{equation}
  N_i \cover{f} N_{i+1}, \quad i \in \mathbb{Z}.
\end{equation}
Then there exists a sequence  $\{z_i\}_{i \in \mathbb{Z}}$ such
that $z_i \in N_i$ and $f(z_i)=z_{i+1}$.
\end{theorem}

The following Lemma plays the crucial role in the construction of
$\rho$ from Theorem~\ref{thm:GH-maps}.
\begin{lemma}
\label{lem:CH-all-cov} Assume the standard conditions.  Let
$$\hat{\alpha}=\hat{\alpha}(A,M)=\max \left(\frac{2M}{c_u -1},
\frac{2M}{1- c_s} \right).$$

Then for any $\alpha >
\hat{\alpha}$, $\lambda_1,\lambda_2 \in [0,1]$ and $z \in
\mathbb{R}^n$ holds that
\begin{equation}
  N(z,\alpha) \cover{A + \lambda_1 h} N((A+\lambda_2h)(z),\alpha).
  \label{eq:GH-all-cov}
\end{equation}
\end{lemma}
\noindent \textbf{Proof:}

Let us fix $z \in \mathbb{R}^n$ and let us define the homotopy
$H:[0,1] \times \overline{B}_u(0,\alpha) \times
\overline{B}_u(0,\alpha)  \to \mathbb{R}^n$ as follows
\begin{equation}
  H_t((x,y))= (A_ux,(1-t)A_sy) + (1-t) \lambda_1 h(z+(r,y))  +
  (A+t\lambda_2 h)(z)
\end{equation}
We have
\begin{eqnarray*}
 H_0(x,y) &=& A(z+(x,y)) + \lambda_1 h(z+(x,y))= (A+\lambda_1
 h)(z+(x,y)) \\
 H_1(x,y) &=& (A+\lambda_2h)(z) + (A_ux,0).
\end{eqnarray*}

For the proof of Lemma~\ref{lem:CH-all-cov} it is enough to show the following conditions for
all $t, \lambda_1,\lambda_2 \in [0,1]$
\begin{eqnarray}
  \| \pi_x (H_t(x,y) - (A+\lambda_2 h)(z) )  \| &>&
  \alpha, \quad (x,y) \in (\partial B_u(0,\alpha)) \times \overline{B}_s (0,\alpha), \label{eq:GH-h-exp-x} \\
   \| \pi_y (H_t(x,y) - (A+\lambda_2 h)(z) )  \| &<&
   \alpha, \quad (x,y) \in \overline{B}_u(0,\alpha) \times \overline{B}_s
   (0,\alpha).  \label{eq:GH-h-contr-y}
\end{eqnarray}
First we establish \eqref{eq:GH-h-exp-x}. We have
\begin{eqnarray*}
    \| \pi_x (H_t((x,y)) - (A+\lambda_2 h)(z) )  \| &=& \\
  \| A_ux + (1-t)\lambda_1\pi_xh(z+(x,y)) + (t-1) \lambda_2 \pi_xh(z) \|
    &\geq& \\
    \| A_ux  \| -  \|h(z+(x,y))\| -   \|h(z)\| &\geq& c_u \alpha - 2M.
\end{eqnarray*}
Hence \eqref{eq:GH-h-exp-x} holds if the following inequality is
satisfied
\begin{equation}
  (c_u - 1) \alpha > 2M.   \label{eq:h-alpha-exp}
\end{equation}

Now we deal with \eqref{eq:GH-h-contr-y}. We have
\begin{eqnarray*}
  \|\pi_y (H_t(x,y) - (A+\lambda_2 h)(z))  \| = \\
  \|(1-t)A_sy + (1-t)\lambda_1 \pi_y h(z+(x,y)) + (t-1) \lambda_2 \pi_y h(z)
  \| \leq
   \\
   \| A_s y\| +   \|h(z+(x,y))\| +  \|h(z)\|   \leq c_s \alpha + 2 M.
\end{eqnarray*}
Hence \eqref{eq:GH-h-contr-y} holds if the following inequality is
satisfied
\begin{equation}
  (1 - c_s ) \alpha > 2M.   \label{eq:h-alpha-contr}
\end{equation}
Hence it is enough take $\hat{\alpha}=\max \left(\frac{2M}{c_u
-1}, \frac{2M}{1- c_s} \right)$.
 \qed

\subsection{The proof of Theorems~\ref{thm:GH-maps} and~\ref{thm:GH-maps-top}}
\label{subsec:thm-gh-maps-proof}

Under assumptions of Theorem~\ref{thm:GH-maps} from Lemma~\ref{lem:homeo} it follows that $g$ is a homeomorphism.
Under assumptions of Theorem~\ref{thm:GH-maps-top}  from Lemma~\ref{lem:onto} it follows that $g$ is a homeomorphism.

 Therefore we can talk of the full orbit of $g$ passing through arbitrary point $z \in \mathbb{R}^n$.

 We define $\sigma:\mathbb{R}^n \to \mathbb{R}^n$  and a multivalued map  $\rho$ from $\mathbb{R}^n$ to subsets of $\mathbb{R}^n$.
 In the case of the proof of Theorem~\ref{thm:GH-maps} $\rho$ we will show that $\rho$ is single valued, i.e. $\rho:\mathbb{R}^n \to \mathbb{R}^n$.
\begin{description}
\item[1] let us fix  $\alpha >\hat{\alpha}$, where $\hat{\alpha}$
is obtained in Lemma~\ref{lem:CH-all-cov},
\item[2] for $z \in \mathbb{R}^n$, from Lemma~\ref{lem:CH-all-cov} with $\lambda_1=1$
and $\lambda_2=0$ we have a bi-infinite chain of covering
relations
\begin{eqnarray}
\dots \cover{g} N(A^{-2}z,\alpha) \cover{g} N(A^{-1}z,\alpha)
\cover{g}  N(z,\alpha) \cover{g} N(Az,\alpha)  \nonumber \\
 \cover{g}
N(A^2z,\alpha) \cover{g}   N(A^3z,\alpha) \cover{g} \dots
\label{eq:GH-chain}
\end{eqnarray}
\item[3.1] in the context of the proof of Theorem~\ref{thm:GH-maps}: from  Theorem~\ref{thm:cov-rel-exists} and Lemma~\ref{lem:bdiff-uniqueness}
it follows that the chain of covering relations  \eqref{eq:GH-chain}  defines
a unique point, which we will denote by
$\rho(z)$, such that
\begin{equation}
   g^k(\rho(z)) \in N(A^k(z),\alpha) \quad k \in \mathbb{Z}.
   \label{eq:GH-def-in}
\end{equation}
\item[3.2] in the context of the proof of Theorem~\ref{thm:GH-maps-top}: from Theorem~\ref{thm:cov-rel-exists} it follows that (\ref{eq:GH-chain})
  defines for each $z \in \mathbb{R}^n$ a non-empty set $\rho(z)$, such that for each $z_1 \in \rho(z)$ holds
  \begin{equation}
   g^k(z_1) \in N(A^k(z),\alpha) \quad k \in \mathbb{Z}.
   \label{eq:GH-def-in-nonunique}
\end{equation}
\item[4] for $z \in \mathbb{R}^n$, from Lemma~\ref{lem:CH-all-cov}
with $\lambda_1=0$ and $\lambda_2=1$ we have a bi-infinite chain
of covering relations
\begin{eqnarray}
\dots \cover{A} N(g^{-2}(z),\alpha) \cover{A} N(g^{-1}(z),\alpha)
\cover{A}  N(z,\alpha) \cover{A} N(g(z),\alpha)  \nonumber \\
 \cover{A}
N(g^2(z),\alpha) \cover{A}   N(g^3(z),\alpha) \cover{A} \dots
\label{eq:GH-inv-chain}
\end{eqnarray}
\item[5] from Theorem~\ref{thm:cov-rel-exists} and the hyperbolicity of $A$
it follows that the chain of covering relations \eqref{eq:GH-inv-chain}
defines a unique point, which we will denote by
$\sigma(z)$, such that
\begin{equation}
   A^k(\sigma(z)) \in N(g^k(z),\alpha) \quad k \in \mathbb{Z}.
   \label{eq:GH-inv-def-in}
\end{equation}
\end{description}

The following lemma shows that in the context of Theorem~\ref{thm:GH-maps} map $\rho$ in fact  does not depend on $\alpha$.
\begin{lemma}
\label{lem:GH-rho-independent-from-alpha}
  Under assumptions of Theorem~\ref{thm:GH-maps}. Assume that $\epsilon < \min(\epsilon_0(A),\epsilon_1(A))$.

  Assume $\hat{\alpha} < \beta$.

  Let $z \in \mathbb{R}^n$. If $z_1$ is such that
  \begin{equation}
     g^k(z_1) \in N(A^kz,\beta), \qquad k \in \mathbb{Z}, \label{eq:gh-gk-beta-in}
  \end{equation}
  then $z_1=\rho(z)$.
\end{lemma}
\textbf{Proof:} Observe that from \eqref{eq:GH-def-in} and
\eqref{eq:gh-gk-beta-in} it follows that
\begin{equation}
  \| g^k(z_1) - g^k(\rho(z))\| \leq \alpha + \beta.
\end{equation}
The assertion follows from Lemma~\ref{lem:bdiff-uniqueness}.
 \qed

The following lemma follows from the hyperbolicity of $A$.
\begin{lemma}
\label{lem:GH-sigma-independent-from-alpha}
 The assumptions as in Theorem~\ref{thm:GH-maps-top}.

  Let $\hat{\alpha} < \beta$.

  Let $z \in \mathbb{R}^n$. If $z_1$ is such that
  \begin{equation}
     A^k(z_1) \in N(g^k(z),\beta), \qquad k \in \mathbb{Z}, \label{eq:gh-Ak-beta-in}
  \end{equation}
  then $z_1=\sigma(z)$.
\end{lemma}

\begin{lemma}
\label{lem:GH-sigma-continous}
 The assumptions as in Theorem~\ref{thm:GH-maps-top}.

 Then $\sigma$ is continuous.
\end{lemma}
\textbf{Proof:}

 Assume that $z_j \to \bar{z}$, we will show that the sequence $\{\sigma(z_j)\}_{j \in
\mathbb{N}}$ is bounded and each  converging subsequence
converges to $\sigma(\bar{z})$.

We can assume that $\|z_j - \bar{z}\| < \alpha$. Then, since
$\|\sigma(z_j) - z_j\| < \alpha$ we obtain
\begin{equation*}
  \| \sigma(z_j) - \bar{z} \| < 2 \alpha.
\end{equation*}
Hence  $\{\sigma(z_j)\}_{j \in \mathbb{N}}$ is bounded.

Now let us take  a convergent subsequence, which we will again
index by $j$, hence $z_j \to \bar{z}$ and $\sigma(z_j) \to w$  for
$j \to \infty$, where $w \in \mathbb{R}^n$.  We will show that $w
= \sigma(\bar{z})$. This implies that $\sigma(z_i) \to \sigma(\bar{z})$.

Let us fix $k \in \mathbb{Z}$. From the continuity of $z \mapsto
g^k(z)$ it follows, that there exists $j_0$ such for $j \geq j_0$
holds
\begin{equation}
   \|g^k(z_j) - g^k(\bar{z})\| < \alpha. \label{eq:Gk-diff}
\end{equation}
Since by the definition of $\sigma$ we have
\begin{equation*}
   A^k(\sigma(z_j)) \in N(g^k(z_j),\alpha)
\end{equation*}
\eqref{eq:Gk-diff} implies that
\begin{equation*}
  \|A^k(\sigma(z_j)) - g^k (\bar{z}) \| \leq 2\alpha.
\end{equation*}
By passing to the limit with $j$  we obtain
\begin{equation}
  \|A^k(w) - g^k (\bar{z})\| \leq 2 \alpha. \label{eq:sigma-cont-limit-k}
\end{equation}
Since \eqref{eq:sigma-cont-limit-k} holds for all $k \in \mathbb{Z}$, then
by Lemma~\ref{lem:GH-sigma-independent-from-alpha}
$w=\sigma(\bar{z})$.
\qed

We continue with the proofs of Theorems~\ref{thm:GH-maps} and~\ref{thm:GH-maps-top}. From the
definition of $\rho$ and $\sigma$ we immediately conclude that  $\sigma \circ g=A \circ \sigma$
and in the context of Theorem~\ref{thm:GH-maps-top} we also havev$\rho \circ A = g \circ \rho$.

We will show that $\sigma(\rho(z))=\{z\}$.

 Let us fix $z  \in \mathbb{R}^n$ and $z_1 \in \rho(z)$, then for any $k \in \mathbb{Z}$ it holds that
\begin{eqnarray*}
  \|g^k(z_1) - A^k(z)\| \leq \alpha, \\
  \|A^k(\sigma(z_1)) - g^k(z_1) \|  \leq \alpha.
\end{eqnarray*}
Hence
\begin{equation}
   \|A^k(\sigma(z_1)) - A^k(z) \|  \leq 2\alpha, \quad k \in
   \mathbb{Z}.
\end{equation}
From the hyperbolicity of $A$ (see also Lemma~\ref{lem:bdiff-uniqueness}) it follows that
$z=\sigma(z_1)$. Therefore we proved
\begin{equation}
\sigma(\rho(z))=\{z\}. \label{eq:sigma-rho}
\end{equation}
Observe that (\ref{eq:sigma-rho}) implies that $\sigma$ is a surjection. This finishes the proof of Theorem~\ref{thm:GH-maps-top}.

From now on we work  under  assumptions of Theorem~\ref{thm:GH-maps} and $\epsilon < \min(\epsilon_0(A),\epsilon_1(A))$.

We will prove that $\rho \circ \sigma=Id$. Let us fix $z \in \mathbb{R}^n$. For all $k \in \mathbb{Z}$ holds
\begin{eqnarray*}
  \|A^k\sigma(z) - g^k(z) \| \leq \alpha, \\
  \|g^k(\rho(\sigma(z))) - A^k\sigma(z)\| \leq \alpha,
\end{eqnarray*}
hence
\begin{equation*}
  \|g^k(\rho(\sigma(z))) - g^k(z)\| \leq 2\alpha.
\end{equation*}
From Lemma~\ref{lem:bdiff-uniqueness} we obtain that $\rho(\sigma(z))=z$.

It remains to show that $\sigma^{-1}=\rho$ is continuous. The proof is virtually the same as the proof of continuity of $\sigma$. The only difference is the use
of Lemma~\ref{lem:GH-rho-independent-from-alpha} in place of Lemma~\ref{lem:GH-sigma-independent-from-alpha}.

\qed

\subsection{Proof of Theorem~\ref{thm:GH-maps-top-no-inj}}
This time we can only consider forward orbits. To define map $\sigma_u$ we proceed as follows.

For any $z \in \mathbb{R}^n$, from Lemma~\ref{lem:CH-all-cov}
with $\lambda_1=0$ and $\lambda_2=1$ we have the following  chain
of covering relations
\begin{eqnarray}
  N(z,\alpha) \cover{A} N(g(z),\alpha)   \cover{A}
N(g^2(z),\alpha) \cover{A}   N(g^3(z),\alpha) \cover{A} \dots
\label{eq:GH-forward-chain}
\end{eqnarray}
From Theorem~\ref{thm:inf-chain} applied to (\ref{eq:GH-forward-chain}) it is easy to show that there exist $z_1=(x_1,y_1) \in \mathbb{R}^u \times \mathbb{R}^s$ such that
\begin{equation*}
  A^k(z_1) \in N(g^k(z),\alpha), \qquad k \in \mathbb{N}.
\end{equation*}

We set
\begin{equation}
  \sigma_u(z)=x_1.
\end{equation}

We need to show first that $\sigma_u(z)$ is well defined. Let $z_2=(x_2,y_2)$ be another point such that
 \begin{equation*}
  A^k(z_2) \in N(g^k(z),\alpha), \qquad k \in \mathbb{N}.
\end{equation*}
Then
\begin{equation}
 \|A_u^k (x_1) - A_u^k(x_2)\| \leq 2\alpha, \quad k \in \mathbb{N}.
\end{equation}
On the other side from our assumptions about $A$ its follows that
\begin{equation}
 \|A_u^k (x_1) - A_u^k(x_2)\| \geq c_u^k  \|x_1 - x_2\|, \quad k \in \mathbb{N}.
\end{equation}
Since $c_u >1$ we conclude that $x_1=x_2$.

From the above reasoning it follows immediately $\sigma_u(z)$ is defined by the following condition
\begin{equation}
  \exists_{\sigma_s(z) \in \mathbb{R}^s} \ \forall k \in \mathbb{N} \quad  A^k(\sigma_u(z),\sigma_s(z)) \in N(g^k(z),\beta), \label{eq:sigma_u-def}
\end{equation}
where $\beta \geq \alpha$.

Let us stress that $\sigma_s(z)$ is not well defined map, there exists many possibilities for $\sigma_s(z)$. However using functional notation $\sigma_s(z)$
will facilitate further discussions.

 To establish the semiconjugacy (\ref{eq:GH-conj-top-no-inj}) observe that from (\ref{eq:sigma_u-def}) we obtain
\begin{equation*}
   \forall k \in \mathbb{N}\setminus \{0\} \quad  A^{k-1}(A(\sigma_u(z),\sigma_s(z)))= A^{k-1}(A_u\sigma_u(z),A_s \sigma_s(z)) \in  N(g^{k-1}(g(z)),\alpha).
\end{equation*}
This implies that
\begin{equation}
  A_u\sigma_u(z)=\sigma_u(g(z)).
\end{equation}

The next step is the continuity of $\sigma_u$.

\begin{lemma}
\label{lem:sigma_u-continous}
 $\sigma_u$ is continuous.
\end{lemma}
\textbf{Proof:}

 Assume that $z_j \to \bar{z}$, we will show that the sequence $\{\sigma_u(z_j)\}_{j \in
\mathbb{N}}$ is bounded and each  converging subsequence
converges to $\sigma_u(\bar{z})$.

We can assume that $\|z_j - \bar{z}\| < \alpha$. Then, since
$\|(\sigma_u(z_j),\sigma_s(z_j)) - z_j\| < \alpha$ we obtain
\begin{eqnarray*}
  \| \sigma_u(z_j) - \pi_x \bar{z} \| < 2 \alpha, \\
   \| \sigma_s(z_j) - \pi_y \bar{z} \| < 2 \alpha.
\end{eqnarray*}
Hence  $\{\sigma_u(z_j),\sigma_s(z_j)\}_{j \in \mathbb{N}}$ is bounded.

Now let us take  a convergent subsequence, which we will again
index by $j$, hence $z_j \to \bar{z}$, $\sigma_u(z_j) \to w$ and $\sigma_s(z_j) \to v$  for
$j \to \infty$, where $w \in \mathbb{R}^u$.  We will show that $w
= \sigma_u(\bar{z})$. This implies that $\sigma_u(z_i) \to \sigma_u(\bar{z})$.

Let us fix $k \in \mathbb{N}$. From the continuity of $z \mapsto
g^k(z)$ it follows, that there exists $j_0$ such for $j \geq j_0$
holds
\begin{equation}
   \|g^k(z_j) - g^k(\bar{z})\| < \alpha. \label{eq:sigmau-Gk-diff}
\end{equation}
Since by the definition of $\sigma_u$ we have
\begin{equation*}
   A^k(\sigma_u(z_j),\sigma_s(z_j)) \in N(g^k(z_j),\alpha)
\end{equation*}
\eqref{eq:sigmau-Gk-diff} implies that
\begin{equation*}
  \|A^k(\sigma_u(z_j),\sigma_s(z_j)) - g^k (\bar{z}) \| \leq 2\alpha.
\end{equation*}
By passing to the limit with $j$  we obtain
\begin{equation}
  \|A^k(w,v) - g^k (\bar{z})\| \leq 2 \alpha. \label{eq:sigmau-cont-limit-k}
\end{equation}
Since \eqref{eq:sigmau-cont-limit-k} holds for all $k \in \mathbb{N}$, then
by (\ref{eq:sigma_u-def}) $w=\sigma_u(\bar{z})$.
\qed

It remains to show the surjectivity of $\sigma_u$.

For this let us set $z=(x_0,0)$ and consider the following chain of covering relations
\begin{equation}
  N(z,\alpha) \cover{g} N(Az,\alpha)  \cover{g}
N(A^2z,\alpha) \cover{g}   N(A^3z,\alpha) \cover{g} \dots
\label{eq:GH-rho-u-chain}
\end{equation}
From Theorem~\ref{thm:inf-chain} applied to (\ref{eq:GH-rho-u-chain}) it follows that there exists $\bar{z}$, such that
\begin{equation*}
  g^k(\bar{z}) \in N(A^k(z),\alpha), \quad k \in \mathbb{N}.
\end{equation*}
Hence
\begin{equation}
  A^k((x_0,0)) \in N(g^k(\bar{z}),2 \alpha), \quad \forall k \in \mathbb{N}. \label{eq:rhu-u-real}
\end{equation}
From (\ref{eq:sigma_u-def}) it follows that
\begin{equation*}
  x_0=\sigma_u(\bar{z}).
\end{equation*}
Since $x_0$ was arbitrary, so $\sigma_u$ is onto.
\qed 

\subsection{From global to local Grobman-Hartman theorem}
\label{sec:local}

The transition from the global to the local version of the Grobman-Hartman theorem is very standard, see for example \cite{Pu,Ze}. We include it here for the sake of completeness sake.

Assume that $\varphi:\mathbb{R}^n \to \mathbb{R}^n$ is a
diffeomorphism satisfying
\begin{equation}
  \varphi(z)=Az + h(z),
\end{equation}
where $A \in \mathbb{R}^{n \times n}$ is a linear hyperbolic
isomorphism and
\begin{equation}
h(0)=0, \quad Dh(0)=0.
\end{equation}

Let us fix $\epsilon>0$. There exists $\delta>0$, such that
\begin{equation}
  \|D h(z)\| < \epsilon, \quad \|z\| \leq \delta.
\end{equation}

Let $t: \mathbb{R}_+ \to \mathbb{R}_+$ be a smooth function such
that
\begin{eqnarray}
  t(r) &=& r, \quad r \leq \delta/2, \\
  t(r) &=& w < \delta, \quad r \geq \delta, \\
  t(r_1) &\leq& t(r_2), \quad r_1 < r_2 \\
  0 &<& t'(r) < 1, \quad r \in [\delta/2,\delta].
\end{eqnarray}

Consider now the function $R: \mathbb{R}^n \to \mathbb{R}^n$ given
by
\begin{equation}
  R(0)=0, \qquad R(z)=\frac{t(\|z\|)z}{\|z\|}, z\neq 0
\end{equation}

It is easy to see that
\begin{eqnarray}
  R(z)&=&z, \quad z \in \overline{B}(0,\delta/2), \\
  R(\mathbb{R}^n) &\subset& \overline{B}(0,w), \\
  \|DR\| &\leq& 1.
\end{eqnarray}

Consider now the following modification of $\varphi$ given by
\begin{equation}
  \hat{\varphi}(z)=Az + h(R(z)).
\end{equation}

It is easy to see that
\begin{eqnarray}
 \hat{\varphi}(z)&=&\varphi(z), \qquad z \in
 \overline{B}(0,\delta_2) \\
 \|h(R(z))\| &\leq& \epsilon \delta , \qquad z \in \mathbb{R}^n, \\
 \|D(h\circ R) (z)\| &\leq & \epsilon, \qquad z \in \mathbb{R}^n.
\end{eqnarray}

It is clear that by taking $\epsilon$ and $\delta$ small enough $h
\circ R $ will satisfy the smallness assumption in
Theorem~\ref{thm:GH-maps} hence we will obtain the local
conjugacy, which is the Grobman-Hartman theorem.

\section{ H\"older regularity of $\rho$}
\label{sec:GH-Holder}

It is know that the conjugating homeomorphism from Theorem~\ref{thm:GH-maps} is H\"older. The geometric proof of this fact is given in the Katok-Hasselblatt book
\cite{KH}. In fact this is a particular case of a more general result about the H\"older regularity of the conjugacy between hyperbolic invariant sets. In \cite{KH} not effort was made to estimate the H\"older exponent in the context of the global Grobman-Hartman theorem.

 Using the functional analysis type approach the H\"older continuity of the conjugating homeomorphism was established by Barreira and Valls \cite{BV}, Bellitskii \cite{B}, Bellitskii and Rayskin \cite{BR} (see \cite{BV} and references given there for other related papers) and apparently the best value of the H\"older exponent was obtained.

 Our goal is to show the H\"older
property for $\rho=\sigma^{-1}$, the map from the conjugacy established in Theorem~\ref{thm:GH-maps}. The main result in this section is
Theorem~\ref{thm:gh-holder}. The same arguments apply also to $\sigma$. We show that we can obtain the same estimate as in \cite{BV,B,BR}.


\begin{lemma}
\label{lem:exp-pos} Let $Q$, $A$, $g$ be as in the proof
of Theorem~\ref{thm:GH-maps}. If $Q(z_1 - z_2) \geq 0$, $z_1 \neq
z_2$. Then $Q(g(z_1) - g(z_2)) > 0$ and
\begin{equation}
 \|\pi_x g(z_1) - \pi_x g(z_2)\| > \theta_u  \|\pi_x z_1 - \pi_x
 z_2\|,
\end{equation}
where $\theta_u=c_u - 2\epsilon_0 > 1$
\end{lemma}
\textbf{Proof:} From the cone condition (Lemma~\ref{lem:cc}) it follows that $Q(g(z_1) -
g(z_2)) > 0$.

Since $Q(z_1 - z_2) \geq 0$, hence
\begin{equation}
  \|\pi_x z_1 - \pi_x z_2\| \geq \|\pi_y z_1 - \pi_y z_2\|.
\end{equation}

We have
\begin{eqnarray*}
 \pi_x g(z_1) - \pi_x g(z_2) =  \int_{0}^1 D \pi_x g(t(z_1-z_2)+z_2)dt
 \cdot (z_1 - z_2) = \\
  A_u  \pi_x (z_1 - z_2) +  \int_{0}^1 \frac{\partial \pi_x h}{\partial x}(t(z_1-z_2)+z_2)dt
 \cdot \pi_x (z_1 - z_2) + \\
  \int_{0}^1 \frac{\partial \pi_x h}{\partial y}(t(z_1-z_2)+z_2)dt \cdot \pi_y (z_1 -
  z_2).
\end{eqnarray*}
Hence for if $Q(z_1 - z_2) \geq 0$ we obtain
\begin{eqnarray*}
 \| \pi_x g(z_1) - \pi_x g(z_2) \| \geq c_u \|\pi_x (z_1 - z_2) \|-
 2 \epsilon \|\pi_x (z_1 - z_2)\|.
\end{eqnarray*}
 \qed

An analogous lemma holds for the inverse map.
\begin{lemma}
\label{lem:exp-neg} Let $Q$, $A$, $g$, $\rho$ be as in the proof
of Theorem~\ref{thm:GH-maps}. If $Q(z_1 - z_2) \leq 0$, $z_1 \neq
z_2$. Then $Q(g^{-1}(z_1) - g^{-1}(z_2)) < 0$ and
\begin{equation}
 \|\pi_y g^{-1}(z_1) - \pi_y g^{-1}(z_2)\| > \theta_s  \|\pi_y z_1 -
 \pi_y
 z_2\|,
\end{equation}
where $\theta_s=\frac{1}{c_s + 2\epsilon} > 1$
\end{lemma}
\textbf{Proof:}
  From the cone condition (Lemma~\ref{lem:cc}) it follows that $Q(g^{-1}(z_1) - g^{-1}(z_2)) < 0$.

Since $Q(z_1 - z_2) \leq 0$, hence
\begin{equation}
  \|\pi_y z_1 - \pi_y z_2\| \geq \|\pi_x z_1 - \pi_x z_2\|.
\end{equation}

We have for any $z_1,z_2$
\begin{eqnarray*}
 \pi_y g(z_1) - \pi_y g(z_2) =  \int_{0}^1 D \pi_y g(t(z_1-z_2)+z_2)dt
 \cdot (z_1 - z_2) = \\
  A_s  \pi_y (z_1 - z_2) +  \int_{0}^1 \frac{\partial \pi_y h}{\partial x}(t(z_1-z_2)+z_2)dt
 \cdot \pi_x (z_1 - z_2) + \\
  \int_{0}^1 \frac{\partial \pi_y h}{\partial y}(t(z_1-z_2)+z_2)dt \cdot \pi_y (z_1 -
  z_2).
\end{eqnarray*}
Hence if $Q(g(z_1)-g(z_2)) \leq 0$, then $Q(z_1 - z_2)<0$ and  we
have
\begin{eqnarray*}
 \| \pi_y g(z_1) - \pi_y g(z_2) \| \leq c_s \|\pi_y (z_1 - z_2) \| +
 2 \epsilon \|\pi_y (z_1 - z_2)\|= \\
 \left( c_s + 2\epsilon \right) \|\pi_y (z_1 - z_2)\|,
\end{eqnarray*}
which after the substitution $z_i \mapsto g^{-1}{z_i}$ gives for
$Q(z_1 - z_2) \leq 0$ the following
\begin{equation}
 \| \pi_y z_1 - \pi_y z_2\| \leq  \left( c_s + 2\epsilon \right)
     \|\pi_y (g^{-1}(z_1) - g^{-1}(z_2))\|.
\end{equation}

 \qed

\begin{lemma}
\label{lem:GH-holder-basic-estm} Let $Q$, $A$, $g$, $\rho$ be as
in the proof of Theorem~\ref{thm:GH-maps}.

Then for any $k \in \mathbb{Z}_+$ holds
\begin{eqnarray}
  \| \rho(z_1) -  \rho(z_2) \| \leq \frac{2 \alpha}{\theta_u^k} +
  \left(\frac{\|A_u\|}{\theta_u}\right)^k \|z_1 - z_2\|, \quad \mbox{if $Q(\rho(z_1)-\rho(z_2)) \geq 0$}, \label{eq:delta-rho} \\
   \| \rho(z_1) -  \rho(z_2) \| \leq \frac{2 \alpha}{\theta_s^k} +
  \left(\frac{\|A_s^{-1}\|}{\theta_s}\right)^k \|z_1 - z_2\|, \quad \mbox{if $Q(\rho(z_1)-\rho(z_2)) \leq 0$}. \label{eq:delta-rho-As}
\end{eqnarray}
\end{lemma}
\textbf{Proof:}
We will consider the case $Q(\rho(z_1)-\rho(z_2))\geq 0$, the case $Q(\rho(z_1)-\rho(z_2))\leq 0$ is analogous, one just need to consider
the inverse maps.

 From  Lemma~\ref{lem:exp-pos} (or Lemma~\ref{lem:exp-neg} in the second case)  applied to $\rho(z_1)$ and $\rho(z_2)$ it
follows that for any $k>0$
\begin{eqnarray*}
  \|g^k(\rho(z_1)) - g^k(\rho(z_2))\| = \|\pi_x g^k(\rho(z_1)) - \pi_x g^k(\rho(z_2))\|
  \geq \\
  \geq \theta_{u}^k \|\pi_x \rho(z_1) - \pi_x \rho(z_2)\| = \theta_u^k \| \rho(z_1) -
  \rho(z_2)\|.
\end{eqnarray*}
Now we derive an upper bound on $  \|g^k(\rho(z_1)) -
g^k(\rho(z_2))\|$. Since $g^k(\rho(z_i)) \in N(A^kz_i, \alpha)$
for $i=1,2$ we obtain
\begin{eqnarray*}
 \|  g^k(\rho(z_1)) - g^k(\rho(z_2)) \| \leq  \|  g^k(\rho(z_1)) - A^kz_1 \| +
   \|A^k z_1 - A^kz_2\| + \\
   \| A^kz_2 -  g^k(\rho(z_2)) \| \leq \alpha + \|A\|^k \|z_1 - z_2\|
   + \alpha = 2\alpha + \|A_u\|^k \|z_1 - z_2\|.
\end{eqnarray*}
By combining the above inequalities  we obtain
\begin{equation}
  \| \rho(z_1) -  \rho(z_2) \| \leq \frac{2 \alpha}{\theta_u^k} +
  \left(\frac{\|A_u\|}{\theta_u}\right)^k \|z_1 - z_2\|. \label{eq:delta-rho-pos}
\end{equation}
\qed


We are now ready to prove the H\"older regularity of $\rho$.
\begin{theorem}
\label{thm:gh-holder} Let $ \gamma =\min\left(\frac{\ln \theta_u}{\ln \|A_u\|},\frac{\ln \theta_s}{\ln \|A_s^{-1}\|}\right)$.
There exists $C>0$, such that any $z_1,z_2 \in \mathbb{R}^n$, $z_1
\neq z_2$  and $\|z_1 - z_2\| < 1$ holds
\begin{equation}
  \frac{\|\rho(z_1) - \rho(z_2)\|}{\|z_1 - z_2\|^\gamma} \leq C,
  \label{eq:GH-holder-reg}
\end{equation}
\end{theorem}
\noindent \textbf{Proof:}
Observe first that $\|A_u\| \geq \theta_u > 1$ and $\|A_s^{-1}\| \geq \theta_s > 1$.

Let us set $\delta_0=1$. Let us denote $\delta=\|z_1 - z_2\|$. For
any $\gamma >0$ and $k \in \mathbb{Z}_+$ from
Lemma~\ref{lem:GH-holder-basic-estm} we have
\begin{equation}
  \frac{\| \rho(z_1) -  \rho(z_2) \|}{\|z_1-z_2\|^\gamma} \leq \frac{2 \alpha}{\theta^k} \delta^{-\gamma} +
  \left(\frac{L}{\theta}\right)^k \delta^{1-\gamma},
\end{equation}
where $(\theta,L)=(\theta_u,\|A_u\|)$ or $(\theta,L)=(\theta_s,\|A^{-1}_s\|)$.

In the sequel we will find $C$ which is good for each case separately, and then we chose the larger $C$.

Observe that \eqref{eq:GH-holder-reg} holds if there exists constants $C_1$ and $C_2$ such that
for each $0< \delta <\delta_0$ there exists $k \in \mathbb{Z}_+$ such that the following
inequalities are satisfied
\begin{eqnarray}
  \frac{2 \alpha}{\theta^k} \delta^{-\gamma} &\leq& C_1,  \label{eq:GH-holder-c1} \\
  \left(\frac{L}{\theta}\right)^k \delta^{1-\gamma} &\leq& C_2. \label{eq:GH-holder-c2}
\end{eqnarray}
We show that we can take
\begin{eqnarray}
C_1=2 \alpha,  \label{eq:C1} \\
 C_2=\frac{L}{\theta}. \label{eq:C2}
\end{eqnarray}

The strategy is as follows: first from \eqref{eq:GH-holder-c1} we
compute $k$ and then we insert it to \eqref{eq:GH-holder-c2},
which will give an inequality, which should hold for any
$0<\delta<\delta_0$, this will produce bound for $\gamma$, $C_1$ and $C_2$.

From \eqref{eq:GH-holder-c1} we obtain
\begin{eqnarray}
  \theta^k &\geq& \frac{2 \alpha \delta^{-\gamma}}{C_1},   \nonumber \\
 k \ln \theta  &\geq& \ln \frac{2\alpha}{C_1} - \gamma \ln \delta.   \label{eq:GH-kC1}
\end{eqnarray}

Taking into account (\ref{eq:C1}) we have
\begin{equation}
  k \ln \theta  \geq  - \gamma \ln \delta.   \label{eq:GH-k}
\end{equation}

We set $k_0=k_0(\delta)= - \frac{\gamma}{\ln \theta} \ln \delta$.
$k_0$ might not belong to $\mathbb{Z}$, but $k_0>0$.  We set
$k=k(\delta)=\lfloor k_0 + 1\rfloor$, where $\lfloor z\rfloor$ is
the integer part of $z$. With this choice of $k$ equation
\eqref{eq:GH-k}  is satisfied. Hence also \eqref{eq:GH-holder-c1} holds.

Now we work on \eqref{eq:GH-holder-c2}. Since
\begin{eqnarray*}
   \left(\frac{L}{\theta}\right)^k \leq
   \left(\frac{L}{\theta}\right)^{k_0+1},
\end{eqnarray*}
then \eqref{eq:GH-holder-c2} is satisfied if the following inequality holds
\begin{eqnarray*}
   \left(\frac{L}{\theta}\right)^{1 - \frac{\gamma}{\ln \theta} \ln \delta} \delta^{1-\gamma} &\leq&
   C_2.
\end{eqnarray*}
By taking the logarithm of both sides of the above inequality we
obtain
\begin{eqnarray*}
 \left(1- \frac{\gamma}{\ln \theta} \ln \delta \right)   \ln \left(\frac{L}{\theta}\right)
   + (1-\gamma) \ln \delta &\leq& \ln C_2.
\end{eqnarray*}
Finally, after an rearrangement of terms arrive at
\begin{eqnarray*}
  \left(1 - \gamma \left(1+ \frac{\ln \frac{L}{\theta}}{\ln \theta}\right) \right)   \ln
  \delta &\leq& \ln C_2  -  \ln \frac{L}{\theta}.
\end{eqnarray*}
The last inequality should  be satisfied for all $\delta \leq
\delta_0=1$. Therefore,  we need the coefficient on lhs by $\ln\delta$
to be nonnegative and the rhs to be nonnegative. It is easy to see
that rhs is nonnegative with $C_2$ given by (\ref{eq:C2}). For the
lhs observe that
\begin{displaymath}
  1+ \frac{\ln \frac{L}{\theta}}{\ln  \theta} = 1 + \frac{\ln L - \ln \theta}{\ln
  \theta}= \frac{\ln L}{\ln \theta}.
\end{displaymath}
Hence we obtain
\begin{eqnarray*}
  1 - \gamma \frac{\ln L}{\ln \theta} \geq 0
\end{eqnarray*}
and finally
\begin{displaymath}
  \gamma \leq \frac{\ln \theta}{\ln L}.
\end{displaymath}
 \qed

\subsection{Comparison with known estimates}
\label{subsec:comparison-with-BV} In \cite[Theorem 1]{BV} (see also \cite{B,BR}) the
following estimate has been given for the H\"older exponent for
the  $\rho$ and $\rho^{-1}$ if the size of the perturbation goes
to $0$ (we use our notation)
\begin{equation}
 \alpha < \alpha_0=\min \left\{-\frac{\ln r(A_s)}{\ln r(A_s^{-1})},
   -\frac{\ln r(A_u^{-1})}{\ln r(A_u)}\right\}, \label{eq:BV-holder-exp}
\end{equation}
where $r(A)$ denotes the spectral radius of the matrix $A$.

Let us consider our estimate of the H\"older exponent from
Theorem~\ref{thm:gh-holder}. In the limit of vanishing
perturbation we obtain (see Lemma~\ref{lem:exp-pos} and
\ref{lem:exp-neg})
\begin{equation}
  \theta_u=c_u, \quad \theta_s=\frac{1}{c_s}.
\end{equation}
Since from assumptions of Theorem~\ref{thm:gh-holder} it follows
that we can assume that
\begin{equation}
  \frac{1}{c_u} = \|A_u^{-1}\|, \qquad c_s = \|A_s\|
\end{equation}
we obtain
\begin{eqnarray*}
  \frac{\ln \theta_u}{\ln \|A_u\|}=  \frac{\ln \frac{1}{\|A_u^{-1}\|}}{\ln
  \|A_u\|}= - \frac{\ln \|A_u^{-1}\|}{\ln \|A_u\|}, \\
  \frac{\ln \theta_s}{\ln \|A_s^{-1}\|}= \frac{\ln \frac{1}{\|A_s\|}}{\ln
  \|A_s^{-1}\|}= -\frac{\ln\|A_s\|}{\ln  \|A_s^{-1}\|}.
\end{eqnarray*}
Therefore our estimate for the H\"older exponent is
\begin{equation}
  \alpha_1<\min \left\{- \frac{\ln \|A_u^{-1}\|}{\ln \|A_u\|}, -\frac{\ln\|A_s\|}{\ln  \|A_s^{-1}\|} \right\}.
\end{equation}
It differs from (\ref{eq:BV-holder-exp}) by the exchange  of the
spectral radius of matrices in (\ref{eq:BV-holder-exp}) by the
norms of matrices.  It is quite obvious by using the adapted norm we can get arbitrary close
to the bound given by (\ref{eq:BV-holder-exp}). For example, if $A_u$ and $A_s$
are diagonalizable over $\mathbb{R}$ if we define the scalar
product so that the eigenvectors are orthogonal, then we obtain
$\|A_{u,s}^{\pm 1}\|=r(A_{u,s}^{\pm 1})$.

To conclude, we claim that we were able to reproduce the H\"older
exponent from \cite{BV,B,BR}.

\section{Grobman-Hartman Theorem for ODEs}
\label{sec:ODE-GH}

Consider an ode
\begin{equation}
  z'=f(z), \quad z \in \mathbb{R}^n, \label{eq:ode}
\end{equation}
such that $f \in C^1$ and $0$ is a hyperbolic fixed point.

It is well know that the Grobman-Hartman theorem is also valid for (\ref{eq:ode}). It can be obtained from Theorem~\ref{thm:GH-maps} for time one map.
 In this section we would like to give a geometric proof, which will not reduce the proof to the map case, but rather we prefer a clean ODE version.

In such approach, the chain of covering relations along the full
orbit will be replaced by an isolating segment along the orbit of
a fixed diameter in the extended  phase space (i.e. $(t,z) \in \mathbb{R} \times \mathbb{R}^n$). The cone conditions
for maps have also its natural analog, we will demand that
\begin{equation}
  \frac{d}{dt} Q(\varphi(t,z_1) - \varphi(t,z_2)) > 0.
\end{equation}

We will consider an ODE
\begin{equation}
  z'=Az + h(z), \quad z \in \mathbb{R}^n.
\end{equation}
We will have the following set of assumptions on $A$ and $h$, which we will refer to as the \emph{ODE-standard conditions}
\begin{itemize}
\item Assume that $A: \mathbb{R}^n \to
\mathbb{R}^n$ is a linear map of the following form
\begin{equation}
  A(x,y)=(A_u x, A_s y)
\end{equation}
where $n=u+s$, $A_u: \mathbb{R}^u \to \mathbb{R}^u$ and
$A_s:\mathbb{R}^s \to \mathbb{R}^s$ are linear  maps such that
\begin{eqnarray}
  (x,A_u x)  &\geq& c_u \|x\|^2, \quad c_u >0, \quad \forall x \in \mathbb{R}^u \\
  (y,A_s y) &\leq&  -c_s \|y\|^2, \quad c_s >0, \quad \forall y \in
  \mathbb{R}^s.
\end{eqnarray}
\item
  Assume that $h: \mathbb{R}^n \to \mathbb{R}^n $ is of class $C^1$ and there exists $M>0$
such that
\begin{equation}
  \|h(x)\| \leq M, \qquad \forall x \in \mathbb{R}^n
\end{equation}
\end{itemize}

Let $\varphi$ be the (local) dynamical system induced by
\begin{equation}
z'=Az+h(z).
\end{equation}

Here is a global version of Grobman-Hartman Theorem for ODEs,
which is similar in spirit to Theorem~\ref{thm:GH-maps}.

\begin{theorem}
 \label{thm:GH-ODEs}
 Assume ODE-standard conditions. Assume additionally that
\begin{equation}
  \|Dh(x)\| \leq \epsilon, \qquad \forall x \in \mathbb{R}^n.
\end{equation}

Under the above assumptions there exists
$\epsilon_0=\epsilon_0(A)>0$, such that if $\epsilon <
\epsilon_0(A)$, then there exists a homeomorphism
$\rho:\mathbb{R}^n \to \mathbb{R}^n$ such that for any $t \in
\mathbb{R}$ holds
\begin{equation}
   \rho (\exp(At)z)  = \varphi(t,\rho(z)).  \label{eq:GH-odes-conj}
\end{equation}
\end{theorem}

\begin{theorem}
 \label{thm:GH-ODEs-top}
  Assume ODE-standard conditions.

 Then there exists a continuous surjective map
$\sigma:\mathbb{R}^n \to \mathbb{R}^n$ such that for any $t \in
\mathbb{R}$ holds
\begin{equation}
    (\exp(At) \sigma(z))  = \sigma(\varphi(t,z)).  \label{eq:GH-odes-semiconj}
\end{equation}
\end{theorem}

In the sequel for $\lambda \in [0,1]$ by
$\varphi^\lambda:\mathbb{R} \times \mathbb{R}^n \to \mathbb{R}^n$
we will denote the dynamical system induced by
\begin{equation}
z'=f^\lambda(z):=Az + \lambda h(z).
\end{equation}

Before the proof of Theorems~\ref{thm:GH-ODEs} and~\ref{thm:GH-ODEs-top} we need first to
develop some technical tools.
 The basic steps and constructions
used in the proof are given in
Section~\ref{subsec:thm-gh-odes-proof}. We invite the reader to
jump first to this section to see the overall picture of the proof
and then consult other more technical sections when necessary.

\subsection{$\varphi^\lambda$ is a global dynamical system}

\begin{lemma}
Assume ODE-standard conditions.

Then for every $(t,z) \in \mathbb{R} \times \mathbb{R}^n$
$\varphi^\lambda(t,z)$ is defined.
\end{lemma}
\textbf{Proof:} Observe that
\begin{equation}
  \|f_\lambda(z)\| \leq \|A\|  \|z\| + M.
\end{equation}
From this using the Gronwall inequality  we obtain the following
estimate
\begin{equation}
  \|z(t)\| \leq \|z(0)\| e^{\|A\| \cdot |t|} +
  \frac{M}{\|A\|} \left(e^{\|A\| \cdot
  |t|}-1 \right).
\end{equation}
This implies that $\varphi^\lambda(t,z)$ is defined.
 \qed

\subsection{Isolating segment}
\label{subsec:iso-seg}

We assume that the reader is familiar with the notion of the
isolating segment for an ode. It has its origin in the Conley
index theory \cite{C} and was developed in papers by Roman
Srzednicki and his coworkers \cite{S1,S2,S3,SW,WZ}.

Roughly speaking, an isolating segment for a (non-autonomous) ode
is the set in the extended phasespace (i.e. $(t,z) \in \mathbb{R} \times \mathbb{R}^n$), whose boundaries are
sections of the vector field. The precise definition can be found
Appendix~\ref{app:isoseg}.

\begin{lemma}
\label{lem:ODE-iso-seg}
Assume ODE-standard conditions.

 There exists $\hat{\alpha}=\max
\left(\frac{2M}{c_u},\frac{2M}{c_s} \right)$, such that for
$\alpha
> \hat{\alpha}$ and for any $\lambda_1,\lambda_2 \in [0,1]$ and
$z_0 \in \mathbb{R}^n$ the set
\begin{equation*}
   N_{\lambda_1}(z_0,\alpha)=\{ (t,(x,y))\ | \ (x-\pi_x \varphi^{\lambda_1}(t,z_0) )^2 \leq \alpha^2,
     \quad (y-\pi_y \varphi^{\lambda_1}(t,z_0) )^2 \leq \alpha^2  \}
\end{equation*}
with
\begin{eqnarray}
   N^-_{\lambda_1}(z_0,\alpha)=\{ (t,(x,y)) \in   N_{\lambda_1}(z_0,\alpha) \ | \
           (x-\pi_x \varphi^{\lambda_1}(t,z_0) )^2=\alpha^2\}, \\
   N^+_{\lambda_1}(z_0,\alpha)=\{ (t,(x,y)) \in   N_{\lambda_1}(z_0,\alpha) \ | \
           (y-\pi_y \varphi^{\lambda_1}(t,z_0) )^2=\alpha^2\}.
\end{eqnarray}
is an isolating segment for $\varphi^{\lambda_2}$.
\end{lemma}
\textbf{Proof:} Let us introduce the following notation
\begin{eqnarray}
  L^-(t,x,y)= (x-\pi_x \varphi^{\lambda_1}(t,z_0) )^2 -  \alpha^2, \\
  L^+(t,x,y)= (y-\pi_y \varphi^{\lambda_1}(t,z_0) )^2 -  \alpha^2.
\end{eqnarray}

 The outside normal vector field to
$N^-_{\lambda_1}(z_0,\alpha)$ is given by $\nabla L^-$. We have
\begin{eqnarray*}
  \frac{\partial L^-}{\partial t}(t,x,y)&=&  -2(x-\pi_x
  \varphi^{\lambda_1}(t,z_0))\cdot  \pi_x f^{\lambda_1}(\varphi^{\lambda_1}(t,z_0))
  )= \\
  & & -2(x-\pi_x  \varphi^{\lambda_1}(t,z_0))\cdot ( A_u \varphi^{\lambda_1}(t,z_0) +
  \lambda_1 \pi_x h(\varphi^{\lambda_1}(t,z_0)) ) \\
   \frac{\partial L^-}{\partial x}(t,x,y)&=&  2(x-\pi_x
  \varphi^{\lambda_1}(t,z_0)), \\
   \frac{\partial L^-}{\partial y}(t,x,y)&=& 0.
\end{eqnarray*}

We verify the exit condition by checking that for $(t,z) \in
N^-_{\lambda_1}(z_0,\alpha)$ holds $\nabla L^-(t,z) \cdot
(1,f^{\lambda_2}(t,z)) >0$.

We have for $(t,(x,y)) \in N^-_{\lambda_1}(z_0,\alpha)$
\begin{eqnarray*}
  \frac{1}{2}\nabla L^-(t,z) \cdot (1,f^{\lambda_2}(t,z))= \\
  -(x-\pi_x  \varphi^{\lambda_1}(t,z_0))\cdot ( A_u \varphi^{\lambda_1}(t,z_0) +
  \lambda_1\pi_x h(\varphi^{\lambda_1}(t,z_0)) ) + \\
  (x-\pi_x  \varphi^{\lambda_1}(t,z_0)) \cdot (A_u x + \lambda_2 \pi_x h(x,y)) = \\
  (x-\pi_x  \varphi^{\lambda_1}(t,z_0)) \cdot (A_u  (x-\pi_x
  \varphi^{\lambda_1}(t,z_0)))+ \\
  (x-\pi_x \varphi^{\lambda_1}(t,z_0))\cdot ( -\lambda_1 \pi_x
h(\varphi^{\lambda_1}(t,z_0)) + \lambda_2 \pi_x h(x,y)) \geq \\
c_u \alpha^2 - 2\alpha M= \alpha (c_u \alpha - 2M).
\end{eqnarray*}
We see that it is enough to take $\hat{\alpha} > \frac{2M}{c_u}$.

For the verification of the  entry condition we will show that for
$(t,z) \in N^+_{\lambda_1}(z_0,\alpha)$ holds $\nabla L^+(t,z)
\cdot (1,f^{\lambda_2}(t,z)) <0$.

 The outside normal vector field to
$N^+_{\lambda_1}(z_0,\alpha)$ is given by $\nabla L^+$. We have
\begin{eqnarray*}
  \frac{\partial L^+}{\partial t}(t,x,y)&=&  -2(y-\pi_y
  \varphi^{\lambda_1}(t,z_0))\cdot  \pi_y f^{\lambda_1}(\varphi^{\lambda_1}(t,z_0))
  )= \\
  & & -2(y-\pi_y  \varphi^{\lambda_1}(t,z_0))\cdot ( A_s \varphi^{\lambda_1}(t,z_0) +
  \lambda_1 \pi_y h(\varphi^{\lambda_1}(t,z_0)) ) \\
    \frac{\partial L^+}{\partial x}(t,x,y)&=& 0, \\
   \frac{\partial L^+}{\partial y}(t,x,y)&=&  2(y-\pi_y  \varphi^{\lambda_1}(t,z_0)).
\end{eqnarray*}

We have for $(t,(x,y)) \in N^+_{\lambda_1}(z_0,\alpha)$
\begin{eqnarray*}
  \frac{1}{2}\nabla L^+(t,z) \cdot (1,f^{\lambda_2}(t,z))= \\
  -(y-\pi_y  \varphi^{\lambda_1}(t,z_0))\cdot ( A_s \varphi^{\lambda_1}(t,z_0) +
  \lambda_1\pi_y h(\varphi^{\lambda_1}(t,z_0)) ) + \\
  (y-\pi_y  \varphi^{\lambda_1}(t,z_0)) \cdot (A_y y + \lambda_2 \pi_y h(x,y)) = \\
  (y-\pi_y  \varphi^{\lambda_1}(t,z_0)) \cdot (A_s  (y-\pi_y \varphi^{\lambda_1}(t,z_0)))+ \\
  (y-\pi_y \varphi^{\lambda_1}(t,z_0))\cdot
  ( -\lambda_1 \pi_y h(\varphi^{\lambda_1}(t,z_0)) + \lambda_2 \pi_y h(x,y)) \leq \\
  -c_s \alpha^2 + 2\alpha M= \alpha (-c_s \alpha + 2M).
\end{eqnarray*}
We see that it is enough to take $\hat{\alpha} > \frac{2M}{c_s}$.

\qed

The following theorem will be obtained using the ideas from the proof of the
Wazewski Rectract Theorem \cite{waz} (see also \cite{C}). We will present the details.
\begin{theorem}
\label{thm:sol-iso-seg}
Assume ODE-standard conditions.
 Let $\alpha > \hat{\alpha}$, where $\hat{\alpha}$ is defined in Lemma~\ref{lem:ODE-iso-seg}.

  Then for any
$\lambda_1, \lambda_2 \in [0,1]$ and $z_0 \in \mathbb{R}^n$, there
exists $z_1 \in \mathbb{R}^n$, such that for all $t \in \mathbb{R}$
holds
\begin{equation}
  \varphi^{\lambda_2}(t,z_1) \in \varphi^{\lambda_1}(t,z_0) +
  B_u(0,\alpha) \times B_s(0,\alpha). \label{eq:in-seg-inf}
\end{equation}
\end{theorem}
\textbf{Proof:}  We will show that for any $T >0$ there exists
$z_T \in z_0+ \overline{B}_u(0,\alpha)\times
\overline{B}_s(0,\alpha)$ such that
\begin{equation}
  \varphi^{\lambda_2}(t,z_T) \in \varphi^{\lambda_1}(t,z_0) +
  \overline{B}_u(0,\alpha) \times \overline{B}_s(0,\alpha), \quad
  t \in [-T,T]. \label{eq:in-seg-T}
\end{equation}
Observe that once (\ref{eq:in-seg-T}) is established by choosing a
convergent subsequence from $z_n \to \bar{z}$ for $n \in
\mathbb{Z}_+$ we obtain an orbit for $\varphi^{\lambda_2}$
satisfying
\begin{equation}
  \varphi^{\lambda_2}(t,z_1) \in \varphi^{\lambda_1}(t,z_0) +
  \overline{B}_u(0,\alpha) \times \overline{B}_s(0,\alpha). \label{eq:in-seg-inf-cl}
\end{equation}.

 From Lemma~\ref{lem:ODE-iso-seg} it follows that $N_{\lambda_1}(z,\alpha)$ is an isolating segments for $\varphi^{\lambda_2}$ for
 any $\lambda_2$.

Let us fix $T>0$. We define  map $h:[0,2T] \times
\overline{B}_u(0,\alpha) \times \overline{B}_s(0,\alpha) \to
\overline{B}_u(0,\alpha) \times \overline{B}_s(0,\alpha) $ as
follows. Let $\tau:N_{\lambda_1}(z_0,\alpha) \to \mathbb{R} \cup
\{\infty\}$ be the exit time function from isolating segment
$N_{\lambda_1}(z_0,\alpha)$ for the process $\varphi^{\lambda_2}$.
From the properties of the isolating segments (see
Appendix~\ref{app:isoseg}) it follows that this function is
continuous.

The map $h(s,\cdot)$ does the following: in the coordinate frame with
moving origin given by $\varphi^{\lambda_1}(s-T,z_0)$ to a point
$z$ we assign $\varphi^{\lambda_2}(s,z)$ if $s$ is smaller than
the exit time, or the exit point (all in the moving coordinate
frame).

The precise definition of $h$ is as follows: let
\begin{equation}
  i(z)=z+\varphi^{\lambda_1}(-T,z_0), \quad \tau_i(z)=\tau(-T,i(z))
\end{equation}
then
\begin{equation}
  h(s,z)=
  \begin{cases}
      \varphi^{\lambda_2}(s,i(z)) -  \varphi^{\lambda_1}(s-T,z_0),      &
          \text{if $s \geq \tau_i(z)$ }, \\
    \varphi^{\lambda_2}(\tau_i(z),i(z))
     - \varphi^{\lambda_1}(\tau_i(z)-T,z_0) & \text{otherwise}.
   \end{cases}
\end{equation}

To  prove (\ref{eq:in-seg-T}) it is enough to show that there
exists $z \in z_0+ \overline{B}_u(0,\alpha)\times \overline{B}_s(0,\alpha) $ such that
\begin{equation}
 \tau(-T,z+\varphi^{\lambda_1}(-T,z_0)) < 2T. \label{eq:exit-time-small}
\end{equation}
We will reason by  contradiction and assume that no such $z$ exists. Since $N_{\lambda_1}(z_0,\alpha)$ is an isolating segment we see that $h$  satisfies the following conditions
\begin{eqnarray}
  h(2T,z) &\in& (\partial B_u(0,\alpha)) \times
  \overline{B}_s(0,\alpha) \quad  \forall z \in \overline{B}_u(0,\alpha)\times
\overline{B}_s(0,\alpha),\\
  h(0,z)&=&z, \qquad \forall z \in \overline{B}_u(0,\alpha) \times
  \overline{B}_s(0,\alpha) \\
  h(s,z)&=&z, \qquad \forall s \in [0,2T], \quad  \forall z \in (\partial B_u(0,\alpha)) \times
  \overline{B}_s(0,\alpha).
\end{eqnarray}
This implies that $h$ is the deformation retraction of
$\overline{B}_u(0,\alpha)\times \overline{B}_s(0,\alpha)$ onto
$(\partial B_u(0,\alpha)) \times
  \overline{B}_s(0,\alpha)$. This is not possible because the homology groups of
  both spaces are different, hence (\ref{eq:exit-time-small}) is true for some $z$.

Hence we obtained (\ref{eq:in-seg-inf-cl}). To have (\ref{eq:in-seg-inf}) for $z_1$ observe that from Lemma~\ref{lem:ODE-iso-seg}
it follows that $(t,\varphi^{\lambda_2}(t,z_1)) \in \inte N_{\lambda_1}(z,\alpha)$ for all $t \in \mathbb{R}$, otherwise it will leave $ N_{\lambda_1}(z,\alpha)$
forward or backward in time. Therefore (\ref{eq:in-seg-inf}) is satisfied.

This finishes the proof.

\qed

\subsection{Cone condition}
\label{subsec:ode-cc}

The cone condition for ODEs is  treated using the methods from
\cite{ZCC} and the cones are defined in terms of a quadratic form.

In this subsection we work under assumptions of Theorem~\ref{thm:GH-ODEs}.

Let $Q(x,y)=(x,x)  - (y,y)$ be a quadratic form on $\mathbb{R}^n$.

By $Q$  we will also denote a matrix, such that $Q(z)=z^tQz$. In
our case $Q=\left[\begin{array}{cc}
  I_u & 0 \\
  0 & -I_s
\end{array}
\right]$,  where $I_u \in \mathbb{R}^{u \times u}$ and $I_s \in
\mathbb{R}^{s \times s}$ are the identity matrices.

\begin{lemma}
\label{lem:ODE-cc} There exists $\epsilon_0=\epsilon_0(A)>0$ such
that if $\epsilon < \epsilon_0$, then there exists $\eta>0$ such
that for $\lambda \in [0,1]$ holds the following \emph{cone
condition}
 \begin{equation}
  \frac{d}{dt} Q(\varphi^\lambda(t,z_1) - \varphi^\lambda(t,z_2)) \geq \pm \eta
  Q(\varphi^\lambda(t,z_1) -  \varphi^\lambda(t,z_2)), \quad \forall z_1,z_2 \in \mathbb{R}^n.
  \label{eq:ODE-cc}
\end{equation}
\end{lemma}
\textbf{Proof:} It is enough to consider (\ref{eq:ODE-cc}) for
$t=0$. We have
\begin{eqnarray*}
 \frac{d}{dt} Q(\varphi^\lambda(t,z_1) - \varphi^\lambda(t,z_2))_{t=0}
 = \\
= (f^\lambda(z_1)-f^\lambda(z_2))^tQ(z_1-z_2) + (z_1-z_2)^t Q
 (f^\lambda(z_1)-f^\lambda(z_2)) = \\
= (z_1-z_2)^t\left( D(z_1,z_2)^tQ + QD(z_1,z_2)
 \right)(z_1-z_2),
\end{eqnarray*}
where
\begin{eqnarray*}
  D(z_1,z_2)=\int_0^1 Df^\lambda(z_2 + t (z_1-z_2))dt=A + \lambda \int_0^1Dh(z_2 + t (z_1-z_2))dt
\end{eqnarray*}
We set
\begin{displaymath}
  C(z_1,z_2)=\int_0^1Dh(z_2 + t (z_1-z_2))dt,
\end{displaymath}
hence
\begin{equation}
 D(z_1,z_2)=A + \lambda C(z_1,z_2), \qquad \|C(z_1,z_2)\| \leq
 \epsilon.
\end{equation}

It is enough to prove that $D^tQ + QD$ is positive definite.
Observe first that $A^tQ + QA$ is positive definite. Indeed, we
have for any $z=(x,y) \in \mathbb{R}^n$
\begin{eqnarray*}
 v^t(A^tQ+QA)v= v^t \cdot \begin{pmatrix}
   A_u^t + A_u & 0 \\
   0 &  -(A_s^t + A_s) \
 \end{pmatrix}  \cdot v = \\
  x^t(A_u^t + A_u)x - y^t(A_s^t +
 A_s)y=2(x,A_ux)-2(y,A_sy) \geq \\
  2c_u x^2 + 2c_s y^2 \geq 2 \min(c_u,c_s) \|v\|^2.
\end{eqnarray*}

Since being a positive definite is an open property we see that
the desired $\eta>0$ and $\epsilon_0 >0$ exist. \qed

\begin{lemma}
\label{lem:ode-bnd-uniqueness} Assume that  $\epsilon <\epsilon_0$
 as in Lemma~\ref{lem:ODE-cc}.  Let $\lambda \in [0,1]$.

Assume that for some $z_1,z_2 \in \mathbb{R}^n$ there exists
$\beta$, such that for all $t \in \mathbb{R}$ holds
\begin{equation}
  \|\varphi^\lambda(t,z_1) - \varphi^\lambda(t,z_2)\| \leq \beta.
\end{equation}
Then $z_1=z_2$.
\end{lemma}
\textbf{Proof:}

Observe that from our assumption it follows that there exists
$\beta_1$, such that
\begin{equation}
 |Q(\varphi^\lambda(t,z_1)-\varphi^\lambda(t,z_2))| \leq \beta_1,
 \qquad \forall t \in \mathbb{R}. \label{eq:absQbnd}
\end{equation}

 We consider two cases: $Q(z_1-z_2)\geq 0$ and
$Q(z_1-z_2)<0$.

Consider first $Q(z_1-z_2)\geq 0$. From Lemma~\ref{lem:ODE-cc} it
follows that for all $t>0$ holds
$Q(\varphi^\lambda(t,z_1)-\varphi^\lambda(t,z_2))> 0$ and for any
$t_0,t >0$ holds
\begin{equation}
 Q(\varphi^\lambda(t+t_0,z_1)-\varphi^\lambda(t+t_0,z_2)) \geq
 \exp(\eta t)
 Q(\varphi^\lambda(t_0,z_1)-\varphi^\lambda(t_0,z_2)).
\end{equation}
This is in a contradiction with (\ref{eq:absQbnd}).

Now we consider case $Q(z_1-z_2)<0$. It is easy to see that
$Q(\varphi^\lambda(t,z_1)-\varphi^\lambda(t,z_2))< 0$ for $t <0$.

From the cone condition (Lemma~\ref{lem:ODE-cc}) it follows that
\begin{equation}
  Q(\varphi^\lambda(t,z_1)-\varphi^\lambda(t,z_2)) <
 \exp(-\eta t) Q(z_1-z_2), \quad t < 0.
\end{equation}
Hence
\begin{equation}
    |Q(\varphi^\lambda(t,z_1)-\varphi^\lambda(t,z_2))| >
 \exp(\eta |t|) |Q(z_1-z_2)|, \qquad t<0.
\end{equation}
This is in a contradiction with (\ref{eq:absQbnd}).

This finishes the proof.
 \qed

\subsection{Proof of Theorems~\ref{thm:GH-ODEs} and ~\ref{thm:GH-ODEs-top}.}
\label{subsec:thm-gh-odes-proof}

The proof follows the pattern of the proof of
Theorems~\ref{thm:GH-maps} and~\ref{thm:GH-maps-top}. Below we will just list the basic steps
of the proof.

 We define $\sigma:\mathbb{R}^n \to \mathbb{R}^n$  and a multivalued map  $\rho$ from $\mathbb{R}^n$ to subsets of $\mathbb{R}^n$.
 In the case of the proof of Theorem~\ref{thm:GH-maps} $\rho$ we will show that $\rho$ is single valued, i.e. $\rho:\mathbb{R}^n \to \mathbb{R}^n$.
\begin{description}
\item[1] let us fix  $\alpha >\hat{\alpha}$, where $\hat{\alpha}$
is obtained in Lemma~\ref{lem:ODE-iso-seg},
\item[2] for $z \in \mathbb{R}^n$, from Lemma~\ref{lem:ODE-iso-seg} with $\lambda_1=1$
and $\lambda_2=0$ we have an isolating segment $N_0(z,\alpha)$ for
$\varphi^1$.
\item[3.1] in the context of the proof of Theorem~\ref{thm:GH-ODEs}: from Theorem~\ref{thm:sol-iso-seg} and Lemma~\ref{lem:ode-bnd-uniqueness} it
follows that $N_0(z,\alpha)$ defines a unique point, which we will
denote by $\rho(z)$, such that
\begin{equation}
   \varphi^1(t,\rho(z)) \in B(\varphi^0(t,z),\alpha) \quad t \in \mathbb{R}.
   \label{eq:GH-odes-def-in}
\end{equation}
\item[3.2] in the context of the proof of Theorem~\ref{thm:GH-ODEs-top}: from Theorem~\ref{thm:sol-iso-seg} it follows that $N_0(z,\alpha)$
  defines for each $z \in \mathbb{R}^n$ a non-empty set $\rho(z)$, such that for each $z_1 \in \rho(z)$ holds
 \begin{equation}
   \varphi^1(t,z_1) \in B(\varphi^0(t,z),\alpha) \quad t \in \mathbb{R}.
   \label{eq:GH-odes-def-in-nonunique}
\end{equation}
\item[4] for $z \in \mathbb{R}^n$, from Lemma~\ref{lem:ODE-iso-seg} with $\lambda_1=0$ and $\lambda_2=1$
we have an isolating segment $N_1(z,\alpha)$ for $\varphi^0$,
\item[5] from Theorem~\ref{thm:sol-iso-seg} and  the hyperbolicity of $A$ it follows
that the isolating segment $N_1(z,\alpha)$ defines a unique point,
which we will denote by $\sigma(z)$, such that
\begin{equation}
   \varphi^0(t,\sigma(z)) \in B(\varphi^1(t,z),\alpha) \quad t \in \mathbb{R}.
   \label{eq:GH-inv-ode-def-in}
\end{equation}
\end{description}
The details of the proof are basically the same as in the proofs of
the map case and are left to the reader.

\section{Appendix. h-set and Covering relations}
\label{app:cov-rel}
 The goal of this section is present
 the notions of the h-set and the covering relation, and to state the theorem
about the existence of point realizing the chain of covering
relations.

\subsection{h-sets and covering relations}
\label{subsec:covrel}

\begin{definition} \cite[Definition 1]{ZGi}
\label{def:covrel} An $h$-set, $N$, is a quadruple \\
$(|N|,u(N),s(N),c_N)$ such that
\begin{itemize}
 \item $|N|$ is a compact subset of $\mathbb{R}^n$
 \item $u(N),s(N) \in \{0,1,2,\dots\}$ are such that $u(N)+s(N)=n$
 \item $c_N:\mathbb{R}^n \to
   \mathbb{R}^n=\mathbb{R}^{u(N)} \times \mathbb{R}^{s(N)}$ is a
   homeomorphism such that
      \begin{displaymath}
        c_N(|N|)=\overline{B_{u(N)}} \times
        \overline{B_{s(N)}}.
      \end{displaymath}
\end{itemize}
We set
\begin{eqnarray*}
   \dim(N) &:=& n,\\
   N_c&:=&\overline{B_{u(N)}} \times \overline{B_{s(N)}}, \\
   N_c^-&:=&\partial B_{u(N)} \times \overline{B_{s(N)}}, \\
   N_c^+&:=&\overline{B_{u(N)}} \times \partial B_{s(N)}, \\
   N^-&:=&c_N^{-1}(N_c^-) , \quad N^+=c_N^{-1}(N_c^+).
\end{eqnarray*}
\end{definition}

Hence an $h$-set, $N$, is a product of two closed balls in some
coordinate system. The numbers $u(N)$ and $s(N)$ are called the
nominally unstable and nominally stable dimensions, respectively.
The subscript $c$ refers to the new coordinates given by
homeomorphism $c_N$. Observe that if $u(N)=0$, then
$N^-=\emptyset$ and if $s(N)=0$, then $N^+=\emptyset$. In the
sequel to make notation less cumbersome we will often drop the
bars in the symbol $|N|$ and we will use $N$ to denote both the
h-sets and its support.

Sometimes we will call $N^-$ \emph{the exit set of N} and $N^+$
\emph{the entry set of $N$}.

\begin{definition}\cite[Definition 6]{ZGi}
\label{def:covw} Assume that $N,M$ are $h$-sets, such that
$u(N)=u(M)=u$ and $s(N)=s(M)=s$. Let $f:N \to \mathbb{R}^n$ be a
continuous map. Let $f_c= c_M \circ f \circ c_N^{-1}: N_c \to
\mathbb{R}^u \times \mathbb{R}^s$. Let $w$ be a nonzero integer.
We say that
\begin{displaymath}
  N\cover{f,w} M
\end{displaymath}
($N$ $f$-covers $M$ with degree $w$) iff the following conditions
are satisfied
\begin{description}
\item[1.] there exists a continuous homotopy $h:[0,1]\times N_c \to \mathbb{R}^u \times \mathbb{R}^s$,
   such that the following conditions hold true
   \begin{eqnarray}
      h_0&=&f_c,  \label{eq:hc1} \\
      h([0,1],N_c^-) \cap M_c &=& \emptyset ,  \label{eq:hc2} \\
      h([0,1],N_c) \cap M_c^+ &=& \emptyset .\label{eq:hc3}
   \end{eqnarray}
\item[2.] If $u >0$, then there exists a  map $A:\mathbb{R}^u \to \mathbb{R}^u$, such that
   \begin{eqnarray}
    h_1(p,q)&=&(A(p),0), \mbox{ for $p \in \overline{B_u}(0,1)$ and $q \in
    \overline{B_s}(0,1)$,}\label{eq:hc4}\\
      A(\partial B_u(0,1)) &\subset & \mathbb{R}^u \setminus
      \overline{B_u}(0,1).  \label{eq:mapaway}
   \end{eqnarray}
  Moreover, we require that
\begin{displaymath}
  \deg(A,\overline {B_u}(0,1),0)=w, \label{eq:deg-A}
\end{displaymath}
\end{description}

We will call condition \eqref{eq:hc2} \emph{the exit condition}
and condition \eqref{eq:hc3} will be called \emph{the entry
condition}.

\end{definition}
Note that in the case $u=0$, if $N \cover{f,w} M$, then $f(N)
\subset \inte M$ and  $w=1$.

\begin{rem}If the map $A$  in condition 2 of Def.~\ref{def:covw} is a linear
map, then  condition (\ref{eq:mapaway}) implies, that
\begin{displaymath}
  \deg(A,\overline {B_u}(0,1),0)=\pm 1.
\end{displaymath}
Hence condition (\ref{eq:deg-A}) is in this situation
automatically fulfilled with $w =\pm 1$.

  In fact, this is the most common situation in the applications of
  covering relations.
\end{rem}

Most of the time we will not interested in the value of $w$ in the
symbol $N \cover{f,w} M$ and we will often drop it and write  $N
\cover{f} M$, instead. Sometimes we may even drop the symbol $f$
and write $N \cover{} M$.

\subsection{Main theorem about chains of covering relations}

\begin{theorem}[Thm. 9]\cite{ZGi}
\label{thm:cov}
 Assume $N_i$, $i=0,\dots,k$, $N_k=N_0$ are
$h$-sets and for each $i=1,\dots,k$ we have
\begin{equation}
  N_{i-1} \cover{f_i,w_i} N_{i} \label{eq:dirgcov}
\end{equation}

 Then there exists a point $x \in \inte N_0$, such that
\begin{eqnarray}
   f_i \circ f_{i-1}\circ \cdots \circ f_1(x) &\in& \inte N_i, \quad i=1,\dots,k \\
  f_k \circ f_{k-1}\circ \cdots \circ f_1(x) &=& x
\end{eqnarray}
\end{theorem}
We point the reader to \cite{ZGi} for the proof.

The following result follows from Theorem \ref{thm:cov}.
\begin{theorem}
\label{thm:inf-chain}
Assume that $I=\mathbb{Z}$ or $I=\mathbb{N}$.
 Let $N_i$, $i\in I$   be h-sets.
Assume that for each $i\in I$ we have
\begin{equation}
  N_{i} \cover{f_{i+1},w_{i+1}} N_{i+1}
\end{equation}

Then there exists a sequence  $\{x_i\}_{i \in I}$, such
that $x_i \in \inte N_i$ and
\begin{eqnarray}
   f_{i+1}(x_{i})=x_{i+1}, \quad \forall i \in I.
\end{eqnarray}
\end{theorem}
\textbf{Proof:}
We will consider the case $I=\mathbb{Z}$, the proof for the other case is almost the same.
 For any $k \in \mathbb{Z}_+$ let us consider a
closed loop of covering relations
\begin{equation*}
   N_{-k} \cover{f_{-k+1}} N_{-k+1} \cover{f_{-k+2}} N_{-k+2}
   \cover{} \dots \cover{f_{k-1}} N_{k-1}  \cover{f_{k}}  N_k
   \cover{A_k} N_{-k},
\end{equation*}
where $A_k$ is some  artificial map such that $ N_k
\cover{A_k}N_{-k}$.

From Theorem~\ref{thm:cov} it follows that   there exists a finite
sequence   $\{x_i^k\}_{i=-k,\dots,k}$    such that
\begin{eqnarray}
       x^k_i &\in& \inte N_i, \\
       f_i(x_{i-1}^k)&=&x^k_i, \quad i=-k+1,\dots,k.
\end{eqnarray}
Since $N_i$ are compact, it is easy to  construct a desired
sequence, by taking suitable subsequences.

\qed

\subsection{Natural structure of h-set} Observe that all the
conditions appearing in the definition of the covering relation
are expressed in 'internal' coordinates $c_{N}$ and $c_{M}$. Also
the homotopy is defined in terms of these coordinates.  This
sometimes makes the matter and the notation look a bit cumbersome.
With this in mind we introduce the notion of a 'natural' structure
on h-set.

\begin{definition}
We will say that $N = \{(x_0,y_0)\} + \overline{B}_u(0,r_1) \times
\overline{B}_s(0,r_1) \subset \mathbb{R}^u \times \mathbb{R}^s$ is
an \emph{$h$-set} with a natural structure given by : \newline
 $u(N)=u$,
$s(N)=s$, $c_{N}(x,y)= \left(\frac{x-x_0}{r_1},\frac{y-y_0}{r_2}
\right)$.

\end{definition}

\section{Appendix. Isolating segments for ODEs}
\label{app:isoseg}

Let us consider the differential equation
\begin{equation}
   \dot{x}=f(t,x)  \label{eq:non-auto-ode}
\end{equation}
where $x \in \mathbb{R}^n$ and $f:\mathbb{R} \times \mathbb{R}^n
\to \mathbb{R}^n$ is $C^1$. Let $x(t_{0},x_{0};\cdot)$ be the
solution of (\ref{eq:non-auto-ode}) such that
$x(t_{0},x_{0};t_{0})=x_{0}$ we put
\begin{equation}
\varphi_{(t_{0},\tau)}(x_{0})=x(t_{0},x_{0};t_{0}+\tau).
\end{equation}
The range of $\tau$ for  which $\varphi_{(t_{0},\tau)}(x_0)$ might
depend on $(t_0,x_0)$. $\varphi$  defines a local flow $\Phi$ on
$\mathbb{R}\times \mathbb{R}^n$ by the formula
\begin{equation}
 \Phi_{t}(\sigma,x)=(\sigma+t,\varphi_{(\sigma,t)}(x)).  \label{eq:semiflow}
\end{equation}
In the sequel we will often call  the first coordinate in the
extended phase space $\mathbb{R} \times \mathbb{R}^n$  {\em the
time}.

We use the following notation: by $\pi_{1}:\mathbb{R}\times
\mathbb{R}^n \rightarrow \mathbb{R}$ and $\pi_{2}:\mathbb{R}
\times \mathbb{R}^n \rightarrow \mathbb{R}^n$ we denote the
projections and for a subset $Z\subset \mathbb{R} \times
\mathbb{R}^n$ and $t\in \mathbb{R}$ we put
\begin{equation}
Z_{t}= \{ x \in \mathbb{R}^n: (t,x)\in Z \}.
\end{equation}

Now we are going to state the definition of  a  isolating segment
for (\ref{eq:non-auto-ode}), which is a modification of the notion
of a periodic isolating segment over $[0,T]$ or $T$-periodic
isolating segment in \cite{S1,S2,S3,SW,WZ}.

\begin{definition}
\label{def:blok} Let $(W,W^{-})\subset \mathbb{R}\times
\mathbb{R}^n$ be a pair of subsets. We call $W$ an {\em  isolating
segment} for (\ref{eq:non-auto-ode}) (or  $\varphi$) if:

\begin{description}
\item[(i)]
$(W,W^-)\cap ([a,b]\times \mathbb{R}^n)$ is a pair of compact sets

\item[(ii)]
for every $\sigma \in \mathbb{R}$, $x\in \partial W_{\sigma}$
there exists $\delta>0$ such that for all $t\in (0,\delta)$
$\varphi_{(\sigma,t)}(x)\not\in W_{\sigma+t}$ or
$\varphi_{(\sigma,t)}(x)\in {\rm int} W_{\sigma+t}$,

\item[(iii)]
\begin{eqnarray*}
& W^{-}= \{ (\sigma,x)\in W: \exists \delta>0\ \forall t\in
(0,\delta)\ \varphi_{(\sigma,t)}(x)
\not\in W_{\sigma+t} \}, \\
&  W^{+}:=\cl(\partial W\setminus W^{-})
\end{eqnarray*}

\item[(iv)]
for all $(\sigma,x)\in W^{+}$ there exists $\delta>0$ such that
$\forall t\in (0,\delta)$ holds
\begin{equation}
\varphi_{(\sigma,-t)}(x) \not\in W_{\sigma-t}
\end{equation}

\item[(v)]
there exists $\eta>0$ such that for all $x \in W^-$ there exists
$t > 0$ such that for all $\tau \in (0,t]$  $\Phi_{\tau}(x)\notin
W$ and $\rho(\Phi_t(x),W) > \eta$
\end{description}

\end{definition}

Roughly speaking, $W^-$ and $W^+$ are sections for
(\ref{eq:non-auto-ode}), through which trajectories leave and
enter the segment $W$, respectively.

\begin{definition}
For  the  isolating segment $W$ we define {\em the exit time
function} $\tau_{W,\varphi}$
\begin{displaymath}
 \tau_{W,\varphi} : W \ni (t_0,x_0) \mapsto \sup \{ t \geq 0: \:
    \forall s \in [0,t] \: (t_0+s,\varphi_{(t_0,s)}(x_0)) \in W \} \in [0,\infty]
\end{displaymath}
\end{definition}
By the Wa\.{z}ewski Retract Theorem \cite{waz} the map
$\tau_{W,\varphi}$ is continuous (compare \cite{C}).


\begin{thebibliography}{ZGi}
\bibitem[A]{A} V. I. Arnol'd, \emph{Supplementary Chapters to the Theory of Ordinary
Differential Equations }. Nauka, Moscow (1978)


\bibitem[BV]{BV} L. Barreira, C. Valls, \emph{H\"older Grobman-Hartman
Linearization}, Discrete and Continuous Dynamical Systems 18
(2007), 187--197


\bibitem[B]{B} G. Belitskii, \emph{On the Grobman-Hartman Theorem in the class $C^\alpha$}, preprint

\bibitem[BR]{BR} G. Belitskii, V. Rayskin, \emph{On the  Grobman-Hartman Theorem in  $\alpha$-H\"older Class For Banach Spaces}, preprint


\bibitem[C]{C} C. C. Conley
{\em Isolated Invariant Sets and the Morse Index.} 1978. CBMS vol.
38, Amer. Math. Soc., Providence

\bibitem[C99]{C99} C. Chicone, \emph{Ordinary Differential Equations with Applications}, New
York: Springer-Verlag, 1999

\bibitem[CS]{CS}  C. Chicone, R. Swanson, \emph{Linearization via the Lie Derivative}, Elect J. Diff. Eqs., Monograph
2, 2000

\bibitem[G1]{G1} D. Grobman, \emph{Homeomorphism of systems of differential
equations}, Dokl. Akad. Nauk SSSR, 128 (159), 880--881

\bibitem[G2]{G2} D. Grobman, \emph{Topological classification of neighborhoods of
a singularity in $n$-space}, Mat. Sb. (N.S.), 56 (98) (1962),
77--94

\bibitem[H1]{H1} P. Hartman, \emph{A lemma in the theory of structural stability
of differential equations}, Proc. Amer. Math. Soc. 11(1960),
610--620

\bibitem[H2]{H2} P. Hartman, \emph{On local homeomorphisms of Euclidean spaces},
Bol. Soc. Mat. Mexicana (2), 5 (1960), 220--241

\bibitem[H3]{H3} P. Hartman, \emph{On the local linearization of differential equations},
Proc. Amer. Math. Soc. 14(1963), 568--573

\bibitem[KH]{KH} A. Katok, B. Hasselblatt, Introduction to the Modern Theory of Dynamical Systems,
Encyclopedia of Mathematics and Its Applications, vo. 54, Cambridge University Press


\bibitem[M]{M} J. Moser, \emph{On a theorem of Anosov}, J. Differential
Equations, 5 (1969), 411--440

\bibitem[Pa]{Pa} J. Palis, \emph{On the local structure of hyperbolic points in Banach spaces},
 An. Acad. Brasil. Ci. 40 (1968), 263--266

\bibitem[PdM]{PdM} J. Palis, W. de Melo \emph{Geometric Theory of Dynamical Systems.
An Introduction}, Springer-Verlag, New York, Heidelberg, Berlin
1982

\bibitem[Pu]{Pu} C. Pugh, \emph{On a theorem of P. Hartman}, Amer.
J. Math., 91 (1969), 363--367


\bibitem[S1]{S1} R. Srzednicki,
{\em Periodic and bounded solutions in blocks for time-periodic
nonautonomuous ordinary differential equations.\/} Nonlin.
Analysis, TMA., 1994, 22, 707--737



\bibitem[S2]{S2} R. Srzednicki,
{\em On detection of chaotic dynamics in ordinary differential
equations\/}, Nonlinear Anal. TMA, 1997, 30 No. 8, 4927-4935.


\bibitem[S3]{S3} R. Srzednicki,
{\em On geometric detection of periodic solutions and chaos,\/}
Proceedings of the Conference ``Nonlinear Analysis and Boundary
Value Problems'', CISM Udine, October 2--6, 1995, CISM Courses and
Lectures Vol. 371, Springer-Verlag, Wien, New York 1996, 197-209.

\bibitem[SW]{SW} R. Srzednicki  and  K. W\'{o}jcik,
{\em A geometric method for detecting chaotic dynamics\/}, J.
Diff. Eq., 1997, 135, 66--82

\bibitem[Wa]{waz} T. Wa\.zewski, \emph{Sur un principe topologique de
l'examen de l'allure asymptotique des int\'egrales des \'equations
diff\'erentielles ordinaires}, Ann. Soc. Polon. Math. 20 (1947),
279--313.

\bibitem[WZ]{WZ} K. W\'ojcik, P. Zgliczy\'nski, \emph{Isolating segments,
fixed point index and symbolic dynamics}, J. Differential
Equations 161 (2000), 245--288.


\bibitem[Ze]{Ze} E. Zehnder, \emph{Lectures on Dynamical Systems. Hamiltonian Vector Fields
  and Symplectic Capacities}, EMS Textbooks in Mathematics, EMS
  2010

 \bibitem[ZGi]{ZGi}P. Zgliczy\'{n}ski and M. Gidea, \emph{Covering relations
for multidimensional dynamical systems}, Journal of Differential
Equations 202/1, 33-58 (2004)

\bibitem[ZCC]{ZCC}P. Zgliczy\'{n}ski, \emph{Covering relations, cone
conditions and the stable manifold theorem}, Journal of
Differential Equations. Volume 246 issue 5, 1774-1819 (2009)

\end{thebibliography}
\end{document}